\documentclass{llncs}
\usepackage{epsfig}
\usepackage{amsmath}
\usepackage{amsfonts}
\usepackage{amssymb}

\newcommand{\R}{{I\!\!R}}

\def\R{{\rm I}\! {\rm R}}

\newtheorem{algorithm}[theorem]{Algorithm}

\begin{document}

\pagestyle{headings}

\title{Modelling approach of a near-far-field model for bubble formation and transport}
\author{J\"urgen Geiser and Paul Mertin}
\institute{Ruhr University of Bochum, \\
The Institute of Theoretical Electrical Engineering, \\
Universit\"atsstrasse 150, D-44801 Bochum, Germany \\
\email{juergen.geiser@ruhr-uni-bochum.de}}
\maketitle

\begin{abstract}

In this paper, we present a model based on a 
near-far-field bubble formation.
We simulate the formation of a gas-bubble in a liquid, e.g.,
water and the transportation of such a gas-bubble in the
liquid.
The modelling approach is based on coupling the 
near-field model, which is done by the Young-Laplace equation,
with the far-field model, which is done with a convection-diffusion equation. 
We decouple the small and large time- and space scales
with respect to each adapted model. Such a decoupling allows
to apply the optimal solvers for each near- or far-field
model.
We discuss the underlying solvers and present the numerical
results for the near-far-field bubble formation and transport model.

\end{abstract}



\section{Introduction}

We are motivated to model bubble formation and transport in liquid,
which are applied in controlled production of gas bubbles
in chemical-, petro-chemical-, plasma- or biomedical-processes,
see \cite{gu2011}, \cite{kumar1970}, \cite{yang2007} and \cite{hayashi2015}.

We consider to decompose the formation process of bubbles, we call it near-field approach,
and the transport process of bubbles, we call it far-field approach.
Such a decomposition allows to separate the large scale-dependencies of the
bubble formation, which has smaller time and space scales as the
transport of the bubbles, which applies larger time and space scales,
see \cite{geiser_2016}.
For such a decomposition, we assume that the bubble is formated from an
orifice in a solid surface and submerged in a liquid (viscous Newtonian liquid),
see \cite{simmons2015}. Therefore, the first process (formation) has to be finalized,
when we start with the second process (transport).
Such a decomposition allows to choose the optimal discretization and solver methods,
i.e., we apply fast ODE-solvers for the near-field model and
level-set methods for the far-field model.

The paper is organized as following:
The modelling problems and their solvers are presented in Section \ref{modell}.
The coupling of the models are discussed in Section \ref{coupling}
The numerical experiments are presented in Section \ref{exp}.
In the contents, that are given in Section \ref{concl}, 
we summarize our results.

\section{Mathematical Model}
\label{modell}

\begin{figure}[ht]
\begin{center}  
\includegraphics[width=10.0cm,angle=-0]{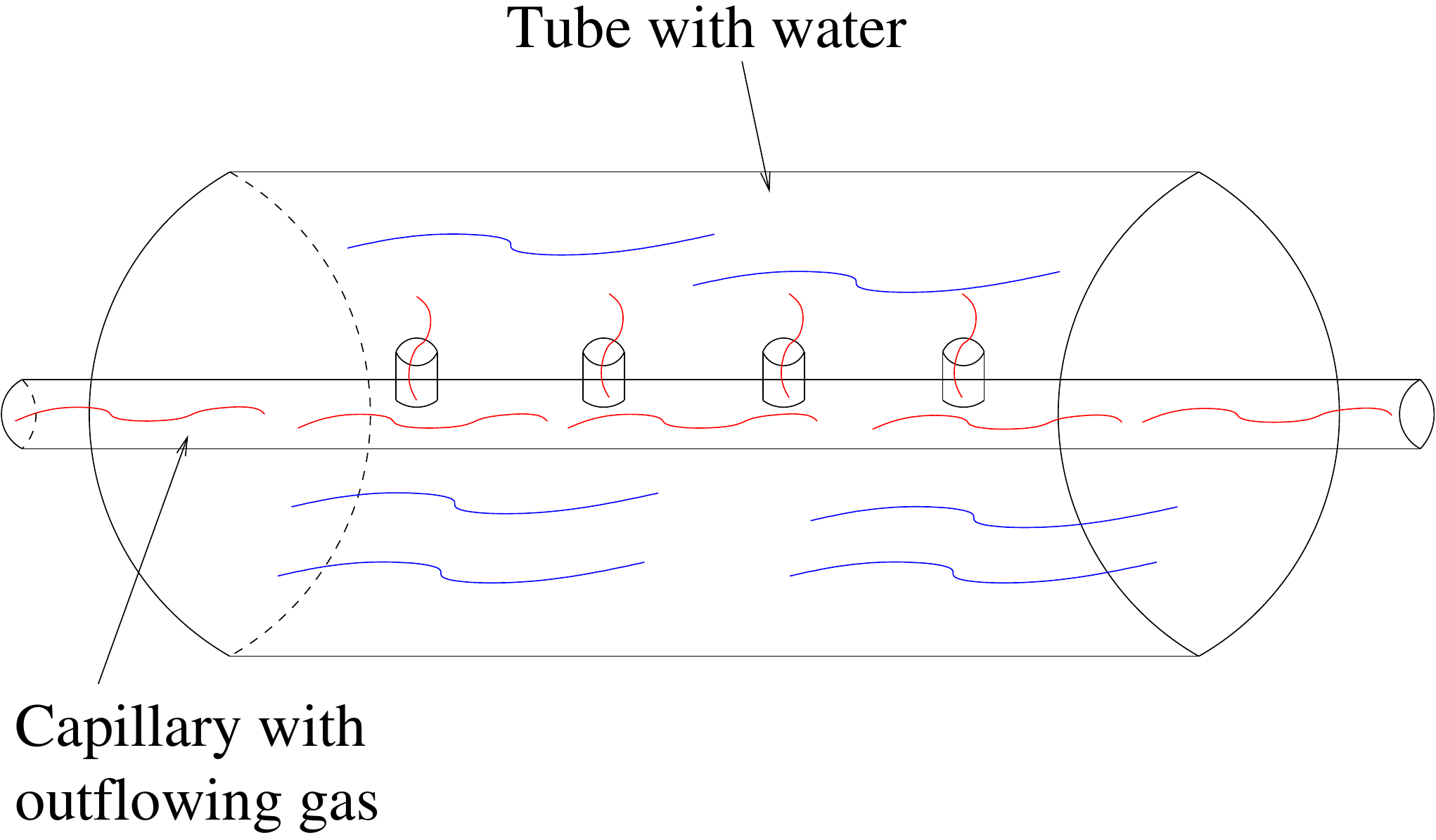}
\end{center}
\caption{\label{experiment} Sketch of the real-life experiment (capillary with gaseous outflow into a tube filled with water).}
\end{figure}

The mathematical model is based on a real-life experiment, where
gas-bubbles are formed in a liquid and are transported after the formation
process, see the plasma-experiment in \cite{hayashi2015}. 
The experiment is given as a thin capillary, where the gas-bubbles are streamed in 
an homogeneous form and transported in a tube, which is filled with
liquid, see the Figure \ref{experiment}.

\begin{figure}[ht]
\begin{center}  
\includegraphics[width=10.0cm,angle=-0]{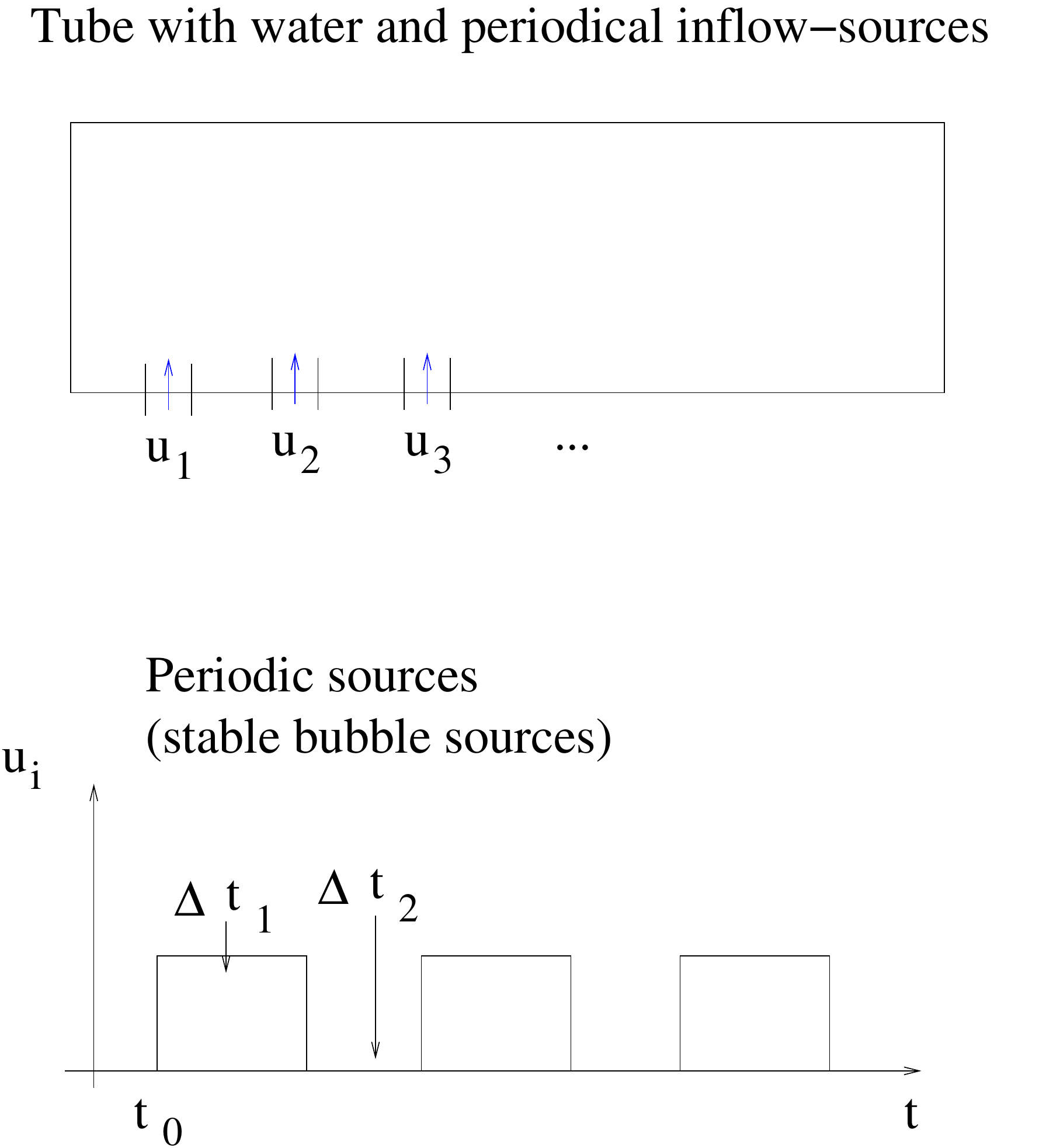}
\end{center}
\caption{\label{bubbles} Periodically inflow of the stable bubble sources.}
\end{figure}

We consider the profile of the tube and deal with the simplified approach of the experiment,
which is given in Figure \ref{bubbles}.

Based on the decoupling of formation and transport, while we assume, that the formation process is not influenced 
by the transport, see \cite{brennen1995}, we deal with two different decoupled models:
\begin{itemize}
\item Near-field approach based on a Young-Laplace equation, see \cite{thorod2005}, where we have a static shape after the formation of the bubbles. 
\item Far-field approach based on a convection-diffusion equation, see \cite{geiser_2016}, where we have a rewriting into a level-set equation, such that we could
transport the static bubble-shapes, see \cite{sethian1996}. 
\end{itemize}

In the following, we discuss the different models.

\subsection{Far-field approach}

The first modelling approach is given with a convection-diffusion equation
in cylindrical coordinate as:
\begin{eqnarray}
&&  \frac{\partial u}{\partial t} = - v \frac{\partial u}{\partial z} +  D_L \frac{\partial^2 u}{\partial z^2} +  \frac{D_t}{r} \frac{\partial}{\partial r} ( r \frac{\partial u}{\partial r}) , \; (r, z, t) \in \Omega \times [0, T] , \\
&&  u(x,z, 0) = u_{near}(x,z) , \; (r, z, t) \in \Omega ,
\end{eqnarray}
where we assume $u_{near}$ is the solution of the bubble-formation in the
near-field and we assume to have Dirichlet-boundary conditions.

Here, we have the benefit and drawbacks of the modelling approach:
\begin{itemize}
\item Benefits:
\begin{itemize}
\item The model is simple and fast to compute.
\item The model also allows to discuss a dynamical shape.
\end{itemize}
\item Drawbacks:
\begin{itemize}
\item The shape of the bubble is not preserved, while we assume a static shape.
\item The influence of the speed of motion in the outer normal direction is not possible. 
\end{itemize}
\end{itemize}

\subsection{Level-Set method}

We apply an improved model, that allows to follow the shapes of the bubble, see \cite{sethian1996}.

The convection-diffusion equation is reformulated in the notation of a level-set equation, which is given as
\begin{align}
& \frac{\partial u}{\partial t}  = - {\bf v} \cdot \nabla u - F_0 |\nabla u| , \; ({\bf x}, t) \in \Omega \times [0, T] , \\
& u({\bf x}, 0) = u_0({\bf x}) , \\
& u(0, t) = u(2, t) = 0.0 , \; ({\bf x}, t ) \in \partial \Omega \times [0, T] ,
\end{align}
where ${\bf v}$ is the convection vector and $F_0$ is the speed of motion in the outer normal direction.
Further $\Omega$ is the computational domain and $T$ is the end time.
The initialisation $u({\bf x}, 0)$ is the results of the near-field computations.

Such equations are wel-known as level-set equations and can be solved like 
convection-diffusion equations, see \cite{sethian1996}.

In the following, we apply the explicit different discretization methods in space, while we
apply the level-set equation with the explicit time-discretisation and apply upwind methods for the advection and outer normal direction term only in the $x$-direction,
the same is also done with the $y$-direction.

We have the following terms:
\begin{align}
& D^-_x u_{i,j} = \frac{u_{i, j} - u_{i-1, j}}{\Delta x} , \\
& D^+_x u_{i,j} = \frac{u_{i+1, j} - u_{i, j}}{\Delta x} , \\
& | D^+_x u_{i,j} | = \left( ( \max(D^-_x u_{i,j}, 0) )^2 + ( \min(D^+_x u_{i,j}, 0) )^2 \right)^{1/2} ,
\end{align}
\begin{align}
 u_i^{n+1}  & = u_i^n -  \Delta t \; v_x \frac{u_{i}^n - u_{i-1}^n}{\Delta x} - \Delta t \; v_y \frac{u_{i,j}^n - u_{i, j-1}^n}{\Delta x} \nonumber \\
& + \Delta t \; F_o | D^+_x u_{i,j} |,
\end{align}
where we assume $v_x, v_y, F_0 \ge 0$.

We discretize the Level-set equation with the explicit time-discretisation and apply upwind methods for the advection and outer normal direction term in the $x$- and $y$-direction
and have the following terms:
\begin{align}
& D^-_x u_{i,j} = \frac{u_{i, j} - u_{i-1, j}}{\Delta x} , \\
& D^+_x u_{i,j} = \frac{u_{i+1, j} - u_{i, j}}{\Delta x} , \\
& D^-_y u_{i,j} = \frac{u_{i, j} - u_{i, j-1}}{\Delta y} , \\
& D^+_y u_{i,j} = \frac{u_{i, j+1} - u_{i, j}}{\Delta y} , \\
& | D^+_x u_{i,j} | = \left( ( \max(D^-_x u_{i,j}, 0) )^2 + ( \min(D^+_x u_{i,j}, 0) )^2 \right. \\
& \left. + ( \max(D^-_y u_{i,j}, 0) )^2 + ( \min(D^+_y u_{i,j}, 0) )^2  \right)^{1/2} .
\end{align}
\begin{align}
 u_i^{n+1}  & = u_i^n -  \Delta t \; v_x \frac{u_{i}^n - u_{i-1}^n}{\Delta x} - \Delta t \; v_y \frac{u_{i,j}^n - u_{i, j-1}^n}{\Delta x} \nonumber \\
& + \Delta t \; F_o | D^+_x u_{i,j} |,
\end{align}
where we assume $v_x, v_y, F_0 \ge 0$.

\begin{remark}
An alternative approach of the shape transport can be done
with the volume of fluid (VOF) method.
Such a method is based on a free-surface modelling technique,
while the method is tracking and locating the free surface, see also \cite{hirt1981}.
\end{remark}

In the following, we discuss the near-field approach.

\subsection{Near field model}

The near-field model is discussed with respect to 
the formation of a drop or bubble, see \cite{thorod2005} and \cite{simmons2015}.

\begin{figure}[ht]
\begin{center}  
\includegraphics[width=5.0cm,angle=-0]{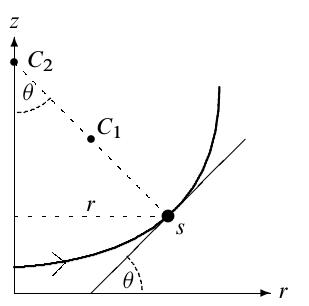}
\end{center}
\caption{\label{near_1} Near field parameters of the bubble shape.}
\end{figure}

The basic modelling idea is based on the so called Young-Laplace equation,
see \cite{finn1986} and deals with the following 
simplied shape of the bubble, see Figure \ref{near_1}.

We deal with the following parameterisation, see \cite{thorod2005}:
\begin{eqnarray}
  && \beta = - \rho g R_t^2 / \sigma , 
\end{eqnarray}
where $\beta$ is the Bond number, $\sigma$ is the surface tension, $\rho$ is the
liquid density, $g$ is the gravity and $R_t$ is the curvature of the drop.

The near-field equations are given as:
\begin{eqnarray}
\label{near_1_1}
  && \frac{d r}{d s} = \cos(\theta) , \\ 
  && \frac{d z}{d s} = \sin(\theta) , \\
\label{near_1_2}
  && \frac{d \theta}{d s} = 2 + \beta z - \frac{\sin(\theta)}{r} , 
\end{eqnarray}
where $s$ is the arc length along the curve
and $\theta$ the angle of elevation for its
slope and $\alpha$ is the mono-layer surface tension.

We have the conditions:
\begin{eqnarray}
  && r=a, \; z=0 , \; \mbox{at} \; s = 0 , \\
  && r = 0, \; \frac{dz}{ds} = \frac{dz}{dr} = 0 , \; \mbox{at} \; s = L ,
\end{eqnarray}
where $L$ is the arc length of the bubble which is a-priori unknown
and is numerically computed via the boundary value problem.
Further, we apply $r = a \rightarrow 0$.

\begin{remark}
The ODE system based on a boundary value problem can be solved with numerical methods, 
e.g., with the MATLAB function bvp4c. Based on the boundary values, we solve the possible curvature of the bubble.
We assume an axial-symmetric bubble and solve the half shape, then we measure the different
diameters of the bubble-ellipse.
\end{remark}

\section{Near-Field Solver: System of ordinary differential equations with boundary conditions}
\label{near-field-solver}

For the near-field, we have to solve a BVP for ODEs,
see \cite{osborne1969} and \cite{liu1993}.

We assume the following nonlinear ODE, give as:
\begin{eqnarray}
  && \frac{d y}{d t} = A(t) y + q(t) , \; t \in [a, b] , \\
  && B_a y(a) + B_b y(b) = d ,
\end{eqnarray}
where ${\bf y}, {\bf q}(t), d \in \R^n$ and $A(t), B_a, B_b \in \R^{n \times}$,
see \cite{liu1993}.

Then, we have the following multiple shooting algorithm, given as

  \begin{algorithm}
\label{algo_1} 

  We have a mesh $a = t_1 < t_2 < \ldots < t_{N+l}= b$, on each mesh interval $[t_i, t_{i+1}]$ with $1 \le i \le N$, we solve
    \begin{enumerate}
   \item A fundamental solution
     \begin{eqnarray}
       Y_i' = A(t) Y_i , \; Y_i(t_i) = F_i , 
     \end{eqnarray}

   \item A particular solution
     \begin{eqnarray}
       {\bf v}_i' = A(t) {\bf v}_i + {\bf q}(t) , \; {\bf v}_i(t_i) = {\bf e}_i , 
     \end{eqnarray}
   \item Then, we find the vector ${\bf s}_i \in \R^n$, $1 \le i \le N$:
     \begin{eqnarray}
       {\bf y}_i' = Y_i(t) {\bf s}_i + {\bf v}_i(t) , \; {\bf v}_i(t_i) = {\bf e}_i , 
     \end{eqnarray}

     \end{enumerate}

  \end{algorithm}

\begin{remark}
The algorithm \ref{algo_1} allows to compute the boundary value problem also
with respect to the Neumann-boundary conditions. Based on the iterative scheme, we also solve the nonlinearity of the equation-system.
\end{remark}

\section{Coupling Near-Field and Far-Field}
\label{coupling}

The modelling assumes, that we could decouple the
near and far-field, while we neglect the coalescence or ruptures
of the bubbles, e.g., in the flow-field, see \cite{leighton1994}.

We assume that in terms of the 
bubble-density function:
\begin{eqnarray}
  && f_b(r, z, x, y, t) = u(x, y, t) \delta((r - R(x,y,t)),(z - Z(x,y,t))),
\end{eqnarray}
where $u$ is the concentration of the bubble and $R$ and $Z$ are obtained with
the bubble-formation equations, while $r$ and $z$ are the cylinder coordinates 
of the density function, that we do not have an influence means
$r \approx R$ and $z \approx Z$ for the formation process.

We discuss the following different coupling ideas:
\begin{itemize}
\item Parameters of the ellipse are computed in the near-field and
initialise the far-field bubble.
\item The near-field computation is directly implemented into the
far-field.
\end{itemize}

\subsection{Decoupled computation of Near- and Far-Field}

The near-field bubble is computed with the ODE's given in
(\ref{near_1_1})-(\ref{near_1_2}).

\begin{figure}[ht]
\begin{center}  
\includegraphics[width=5.0cm,angle=-0]{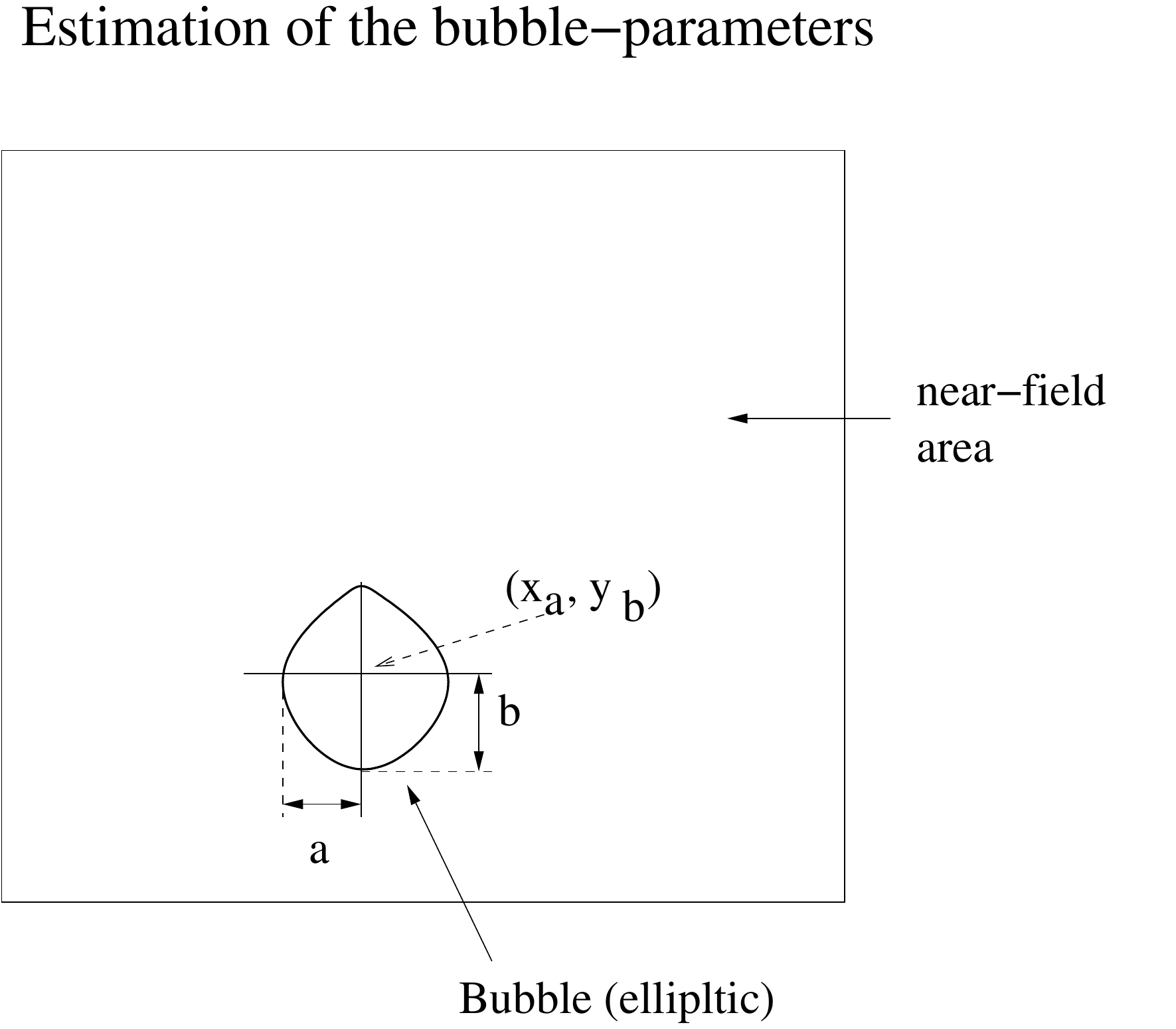}
\end{center}
\caption{\label{ellipse_1} Final bubble based on the near-field computation and estimation of the bubble-parameters (we assume an elliptic curve).}
\end{figure}

We estimate the characteristic parameters of the ellipse
in the Figure \ref{ellipse_1}.

Based on the estimation of the elliptic-parameters, we obtain the
curvature of the ellipse:
\begin{eqnarray}
  && \frac{(x - x_a)^2}{a^2} + \frac{(y - y_b)^2}{b^2} = 1 .
\end{eqnarray}

The ellipse is the curvature of the far-field, which is 
computed with the level-set method.

\begin{remark}
The transformation of the elliptic parameters of the
near-field model allows to simplify the construction of the
shape in the far-field. We only apply the ellipses in the
far-field transport modell.
\end{remark}

\subsection{Coupled computation of Near- and Far-Field}

Based on the electric-field on smaller time-scale, we have to update
the near-field and far-field coupling in a numerical cycle,

The near-field bubble is computed with the ODE's and with the
underlying influence of the E-field.
Further the far-field is computed by the transport equation with
level-set methods. We couple the coupling via an interpolation
between the mesh-free space of the near-field and the grid-space of the
far-field with the give diameters of the bubbles and the trajectories
of the transport field

The oscillation of the spherical bubble, see \cite{leighton1994}, \cite{sommers2012} and \cite{bellini1997}, is given as:
\begin{eqnarray}
&& r = r_0 + r_{\epsilon}, \\
&& r_{\epsilon} = - r_{\epsilon_0} \exp(i \; \omega_0 t),
\end{eqnarray}
where $r_0$ is the mean radius, $\omega_0$ is the resonance frequency
and $r_{\epsilon_0}$ amplitude of maximal oscillations of the bubble radius.

The characteristic frequency of the breathing mode with respect
to the electric stress (pressure) is given as:
\begin{eqnarray}
&& \omega_0 = \frac{1}{2 \pi r_0} \sqrt{\frac{3 \; k \; p}{\rho} - \frac{2 \sigma}{\rho \; r_0}},
\end{eqnarray}
where $\sigma$ is the gas-liquid surface tension, $\rho$ is the liquid density,
$k$ is the polytropic exponent ($k = 1.4$ for air). $p$ is the pressure, where $p = |p_0 - p_E|$ with $p_0$ is the hydrostatic pressure in the fluid far from the bubble and $p_E$ is the electric stress (pressure) given of the electric field
with
\begin{eqnarray}
&& p_E = \frac{9}{8} \epsilon |{\bf E}_0|^2 \; \sin^2,
\end{eqnarray}
the pressure maximum is given at $\theta = \pi/2$, $\epsilon$ is the peremitivity
and ${\bf E}_0$ is the uniform electric field.

\begin{figure}[ht]
\begin{center}  
\includegraphics[width=10.0cm,angle=-0]{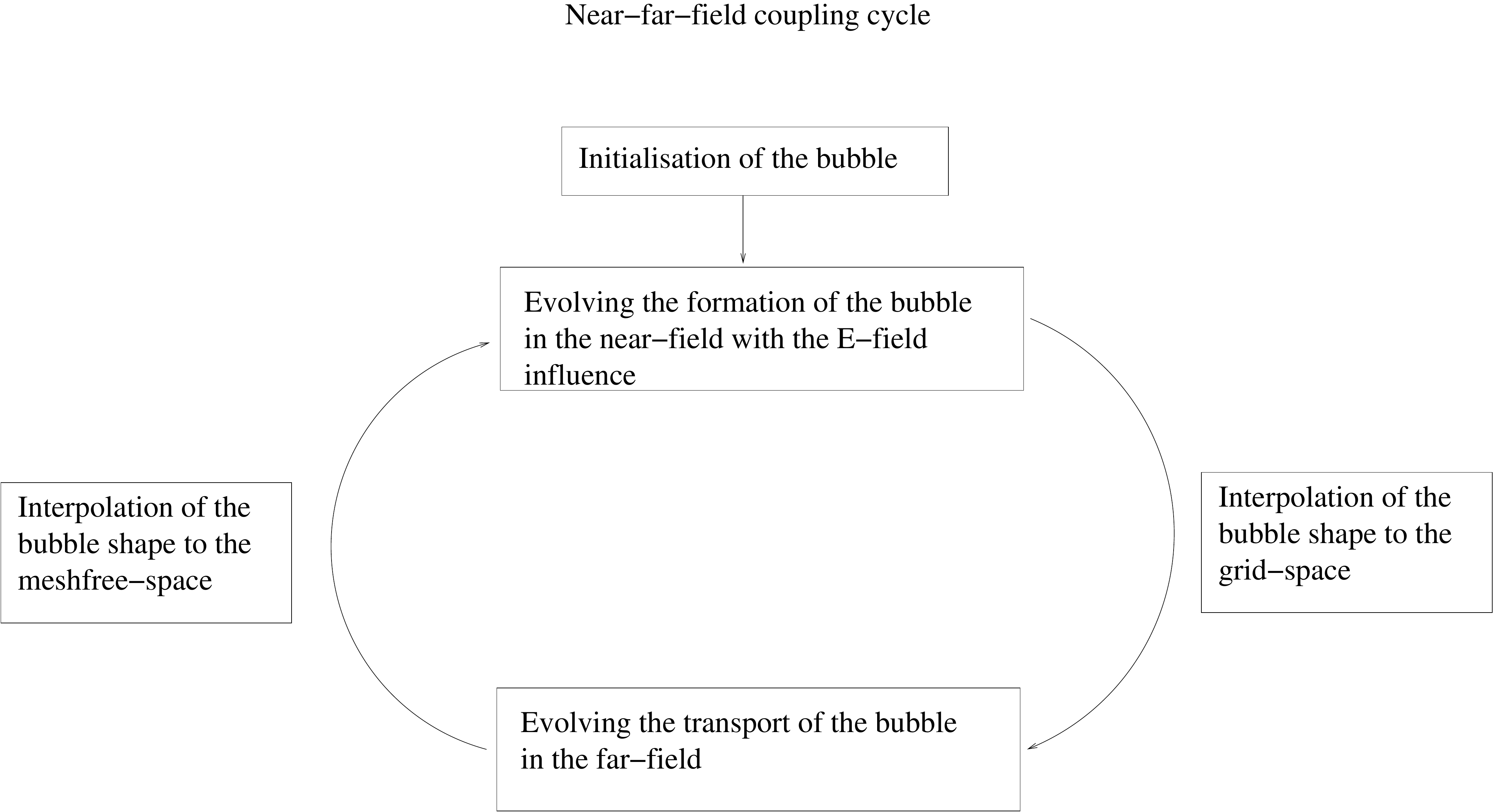}
\end{center}
\caption{\label{cycle} Numerical cycle of the near-far-field coupling.}
\end{figure}

The numerical cycle is given in Figure \ref{cycle}.

\begin{figure}[ht]
\begin{center}  
\includegraphics[width=10.0cm,angle=-0]{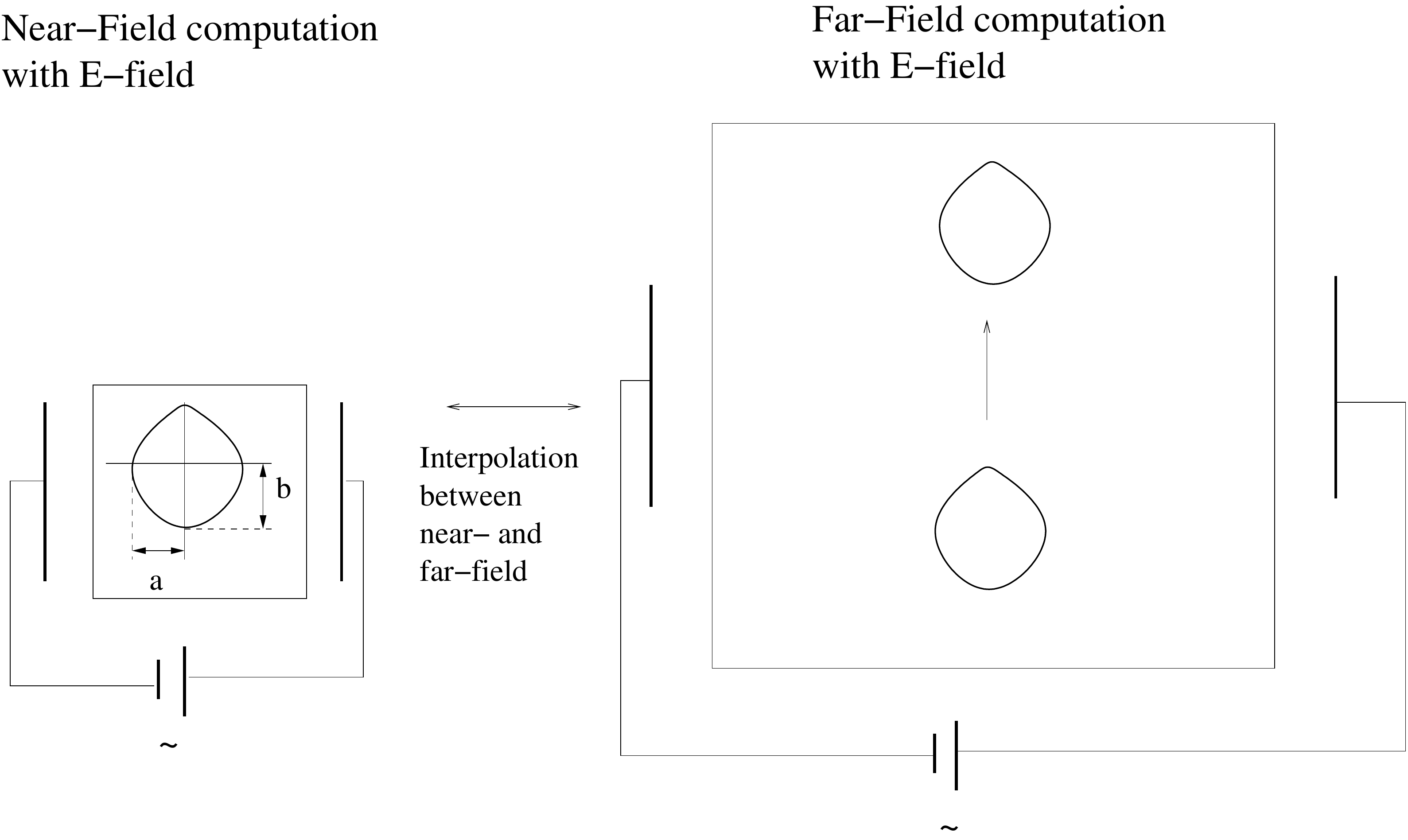}
\end{center}
\caption{\label{cycle_2} Near field and far field coupling in the E-field.}
\end{figure}

In the following, we discuss the detailed coupling of the
near-field and far-field, see Figure \ref{cycle_2}.

\begin{remark}
The accuracy of the cycle is give with respect to the step-size of the
far-field.
\end{remark}

\section{Numerical Experiments}
\label{exp}

In the following, we apply the different numerical experiments
based on the bubble-formation and the bubble-transport.
While the bubble-formation is based on the model with Young-Laplace equation,
the bubble-transport model is based on the
level-set equations.

\subsection{Bubble-Formation: Experiment 1}

The near-field equations are given as:
\begin{eqnarray}
  && \frac{d r}{d s} = \cos(\theta) , \\ 
  && \frac{d z}{d s} = \sin(\theta) , \\
  && \frac{d \theta}{d s} = - \frac{\sin(\theta)}{r} + \frac{\Delta p}{\alpha} , 
\end{eqnarray}
or in vectorial form:
\begin{eqnarray}
  && \frac{d y}{d t} = \tilde{A}(y(t)) + \tilde{q}(t) , \; t \in [t_i, t_{i+1}] ,
\end{eqnarray}
where $\tilde{A}(y(t)) = \left( \begin{array}{c} \cos(\theta) \\ \sin(\theta) \\  - \frac{\sin(\theta)}{r} \end{array} \right) $
and  $\tilde{q}(t) = \left( \begin{array}{c} 0 \\ 0 \\  \frac{\Delta p}{\alpha} \end{array} \right) $.
Further, $y(t) = (r(t), z(t), \theta(t))^t$ and $t$ is the arc length along the curve
and $r$ is the radius, $z$ is the vertical distance and $\theta$ the angle of elevation for its
slope and $\alpha$ is the mono-layer surface tension.

We have to linearise the equations, which are given as:
\begin{eqnarray}
  && \frac{d y}{d t} = \tilde{A}(y(t_i)) + \frac{\partial \tilde{A}}{\partial y}|_{y(t_i)} (y(t) - y(t_i)) + \tilde{q}(t) , \; t \in [t_i, t_{i+1}] , \\
  && \frac{d y}{d t} = A(t) y(t) + q(t) , \; t \in [t_i, t_{i+1}] , \\
  && B_a y(a) + B_b y(b) = d ,
\end{eqnarray}
where $A(t) = \frac{\partial \tilde{A}}{\partial y}|_{y(t_i)} y(t) $, $q(t) = \tilde{A}(y(t_i)) - \frac{\partial \tilde{A}}{\partial y}|_{y(t_i)} y(t_i) + \tilde{q}(t)$. Further, we have $y = (r, z, \theta)^t$ and
\begin{eqnarray}
  && \frac{\partial A}{\partial y}|_{y(t_i)} = \left( \begin{array}{ccc}
0 & 0 & - \sin(\theta(t_i)) \\
0 & 0 & \cos(\theta(t_i)) \\
- \frac{\sin(\theta(t_i))}{r^2(t_i)} & 0 & - \frac{\cos(\theta(t_i))}{r(t_i)} 
  \end{array} \right) .
\end{eqnarray}

We have the conditions:
\begin{eqnarray}
  && r=a, \; z=0 , \; \theta = \pi/2 , \; \mbox{at} \; t = 0 , \\
  && r = 0 , \; \frac{dz}{ds} = 0 , \; \theta = 2 (\pi/2) , \mbox{at} \; t = L  ,
\end{eqnarray}
where $a$ is the radius of the bubble, $L$ is the arc length of the bubble which is a-priori unknown, where we start with $L = L_0 = \frac{2 \pi a^2}{4}$ and we go on with $L = L_1 \le L_2 \le \ldots \le L_{stop}$.

We have the following linear equation system:
\begin{eqnarray}
  && \left[ \begin{array}{llllll}
    S_1 & R_1 &           \\
        & S_2 & R_2 &     \\
    &     & \ddots & \ddots &  \\
    &     &        &  S_N &  R_N \\
  B_a &   &        &   B_{b-1}   &  B_b  
    \end{array} \right]
   \left[ \begin{array}{l}
    y_1 \\
    y_2 \\
    \vdots \\
    y_N \\
    y_{N+1}  
     \end{array} \right] =
    \left[ \begin{array}{l}
    q_1 \\
    q_2  \\
    \vdots \\
    q_N \\
    d  
    \end{array} \right]
\end{eqnarray}4
where $S_i = R_i = - h_i^{-1} I - \frac{1}{2} A(t_{i + 1/2})$, $q_{i} = q(t_{i+1/2}$
with $h_i = t_{i+1} - t_i$, $t_{i + 1/2} = t_i + \frac{1}{2} h_i $.

Further, we have the matrices and vectors:
\begin{eqnarray}
  && B_a = \left( \begin{array}{ccc}
    0 & 0 & 0 \\
   0 & 0 & 0 \\
   0 & 0 & 0 
  \end{array} \right) , \\
    && B_{b-1}  = \left( \begin{array}{ccc}
    0 & 0 & 0 \\
   0 & 1 & 0 \\
   0 & 0 & 0 
  \end{array} \right) , \\
     && B_{b}  = \left( \begin{array}{ccc}
   1 & 0 & 0 \\
   0 & 1 & 0 \\
   0 & 0 & 1 
  \end{array} \right) , \\
     && d  = \left( \begin{array}{ccc}
   L \\
   0 \\
   0
  \end{array} \right) .
\end{eqnarray}

\begin{itemize}
  \item Test example 1: \\ 
We apply the following test-example with the
a quarter of a circle, means $L = 2 \pi a^2/ 4$, where $a = 1, 2, 3, 3.5$
and we assume $N = L / \Delta t$, where we assume $\Delta t = 2 \pi \; 0.001$ and $\frac{\Delta p}{\alpha} = 0$.
\item Test example 2: \\
We have $N = L / \Delta t$, where we assume $\Delta t = 0.001$ and $L = 2$
and $\frac{\Delta p}{\alpha} = 0, 0.8, 1.4$.
\end{itemize}

\begin{figure}[ht]
\begin{center}  
\includegraphics[width=5.0cm,angle=-0]{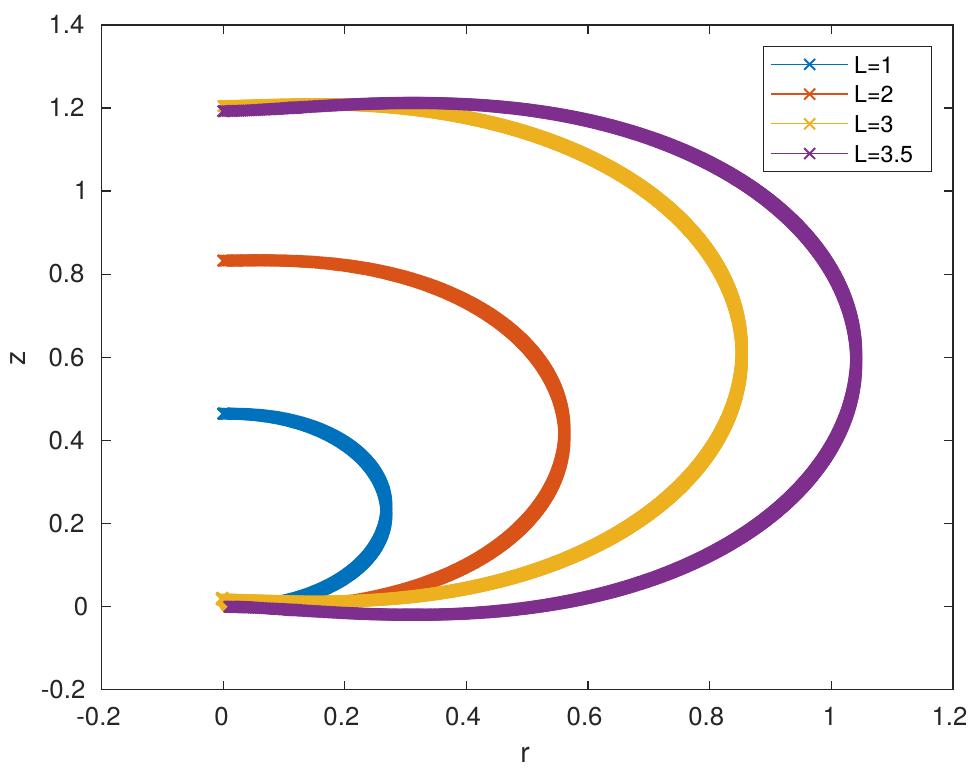}
\includegraphics[width=5.0cm,angle=-0]{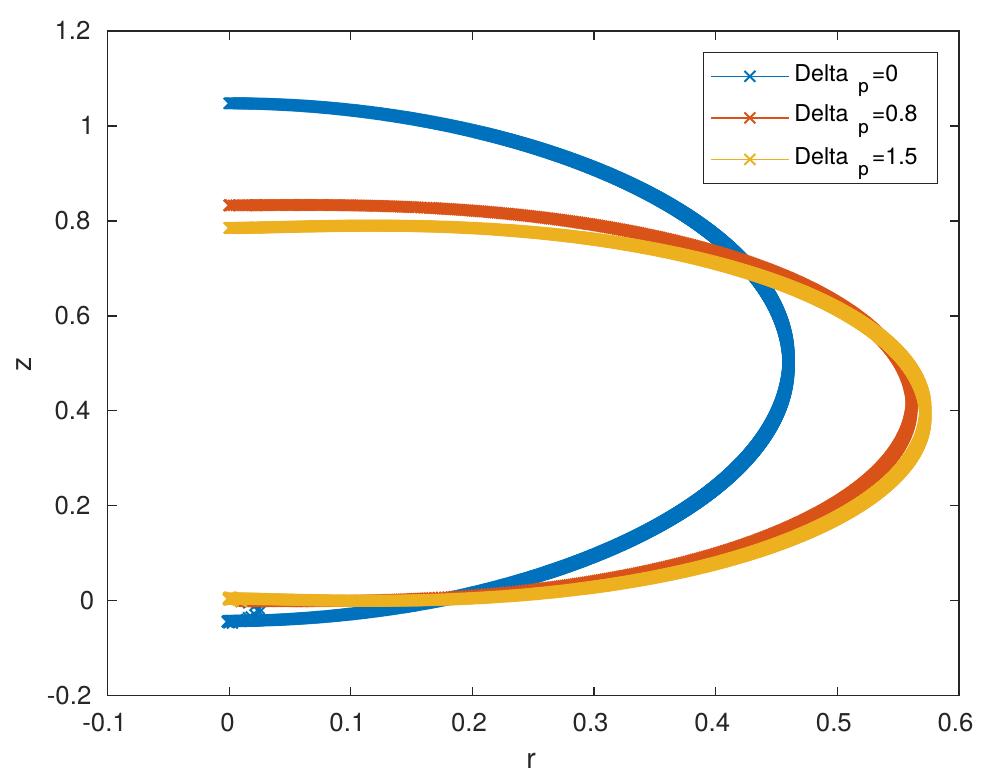}
\end{center}
\caption{\label{bubble_1} Bubble-formation, left figure with different length $L$ and $\frac{\Delta p}{\alpha} = 0$, while $L = 2$ is exact solution of the $1/4$ circle and right figure with different right-hand side-parameters $\frac{\Delta p}{\alpha}$ and for $L = 2$.}
\end{figure}

The numerical results of bubble-formation is given in Figure \ref{bubble_1}.

\begin{remark}
The Young-Laplace equation allows to formulate the bubble-formation such that we could obtain the radii of the different bubbles based on the various pressure parameters.
We also compare the results with respect to the computations in the Literature, see \cite{simmons2015}.
\end{remark}

\subsection{Bubble-Formation: Experiment 2}

In the following, we couple the near-field and far-field computations.

\begin{figure}[ht]
\begin{center}  
\includegraphics[width=8.0cm,angle=-0]{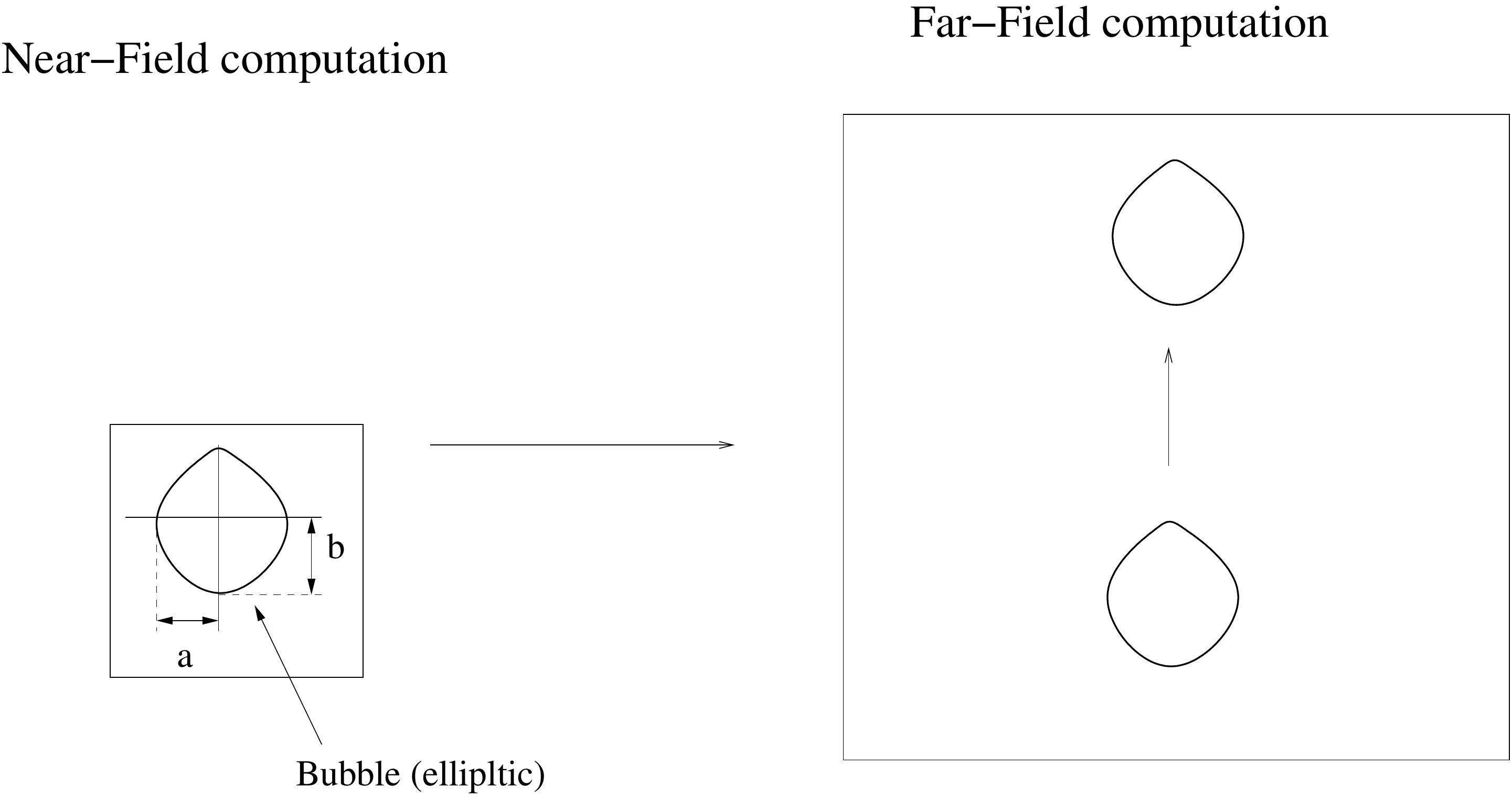}
\end{center}
\caption{\label{near-far_1} Left figure: Bubble-formation with the ODEs and right figure: Bubble-transport with the PDEs (level-set equations).}
\end{figure}

We have the following setting, see Figure \ref{near-far_1}.

We apply the following parameters:
\begin{itemize}
\item Input-parameters of the near-field bubbles:
\begin{itemize}
\item $r_1(0)=0.01, \; z_1(0)=0, \; r_1(L) =0, \; L = 1, \; \frac{\Delta p_1}{\alpha_1}= 0.8$, 
\item $r_2(0)=0.01, \; z_2(0)=0, \; r_2(L) =0, \; L = 2, \; \frac{\Delta p_2}{\alpha_2}= 0.8$.
\end{itemize}
\item We compute the bubbles based on the near-field code and obtain the ellipse-diameters $a_{bubble}, \; b_{bubble}$.
\item We initialise the two ellipses for the far-field computations given as:
\begin{itemize}
\item  $(x-20)^2 + ((y-100)*a_{bubble, 1}/b_{bubble, 1})^2 - a_{bubble, 1}^2 $,
\item  $(x-70)^2 + ((y-90)*a_{bubble, 2}/b_{bubble, 2})^2 - a_{bubble, 2}^2 $.
\end{itemize}
\end{itemize}

\begin{figure}[ht]
\begin{center}  
\includegraphics[width=5.0cm,angle=-0]{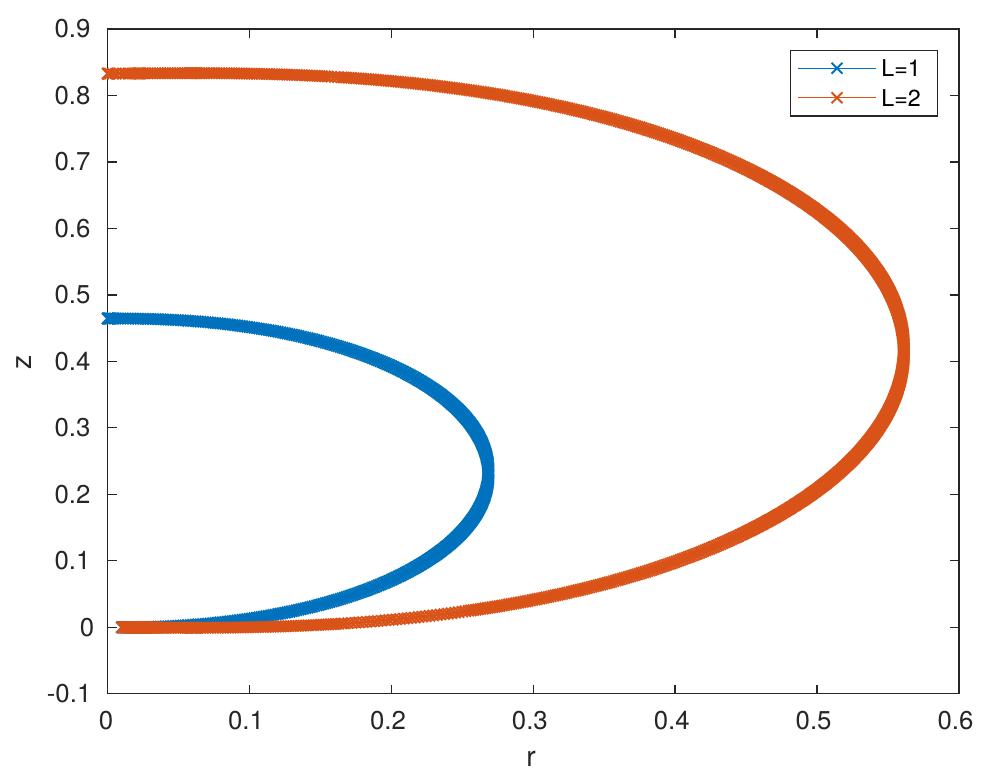}
\end{center}
\caption{\label{bubble_near_2_bubb} The computation of the bubble formation (near-field) for two bubbles.}
\end{figure}

The numerical results of near-field bubble-formation is given in Figure \ref{bubble_near_2_bubb}.

\begin{figure}[ht]
\begin{center}  
\includegraphics[width=5.0cm,angle=-0]{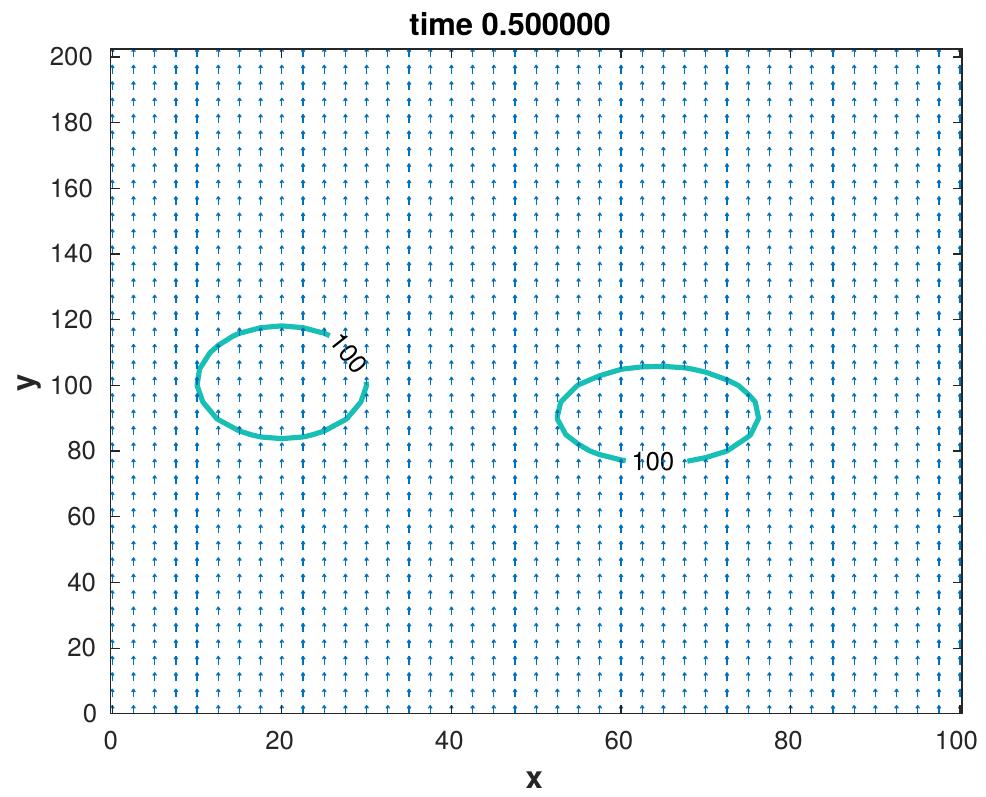}
\includegraphics[width=5.0cm,angle=-0]{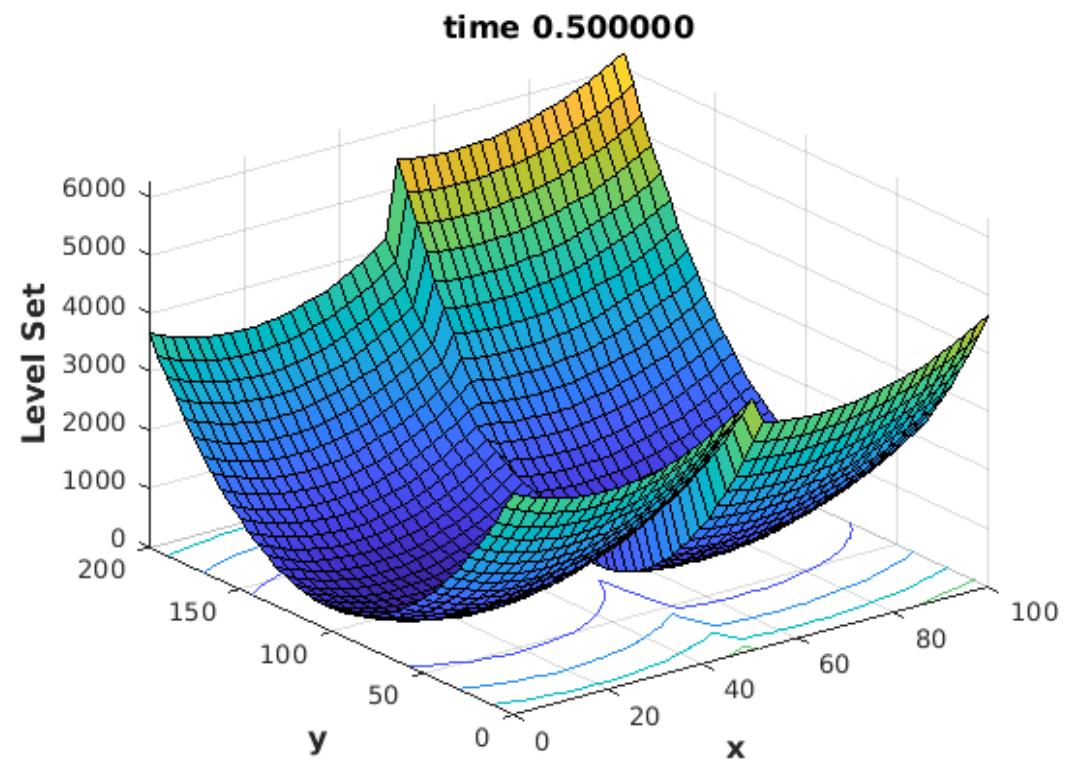} \\
\includegraphics[width=5.0cm,angle=-0]{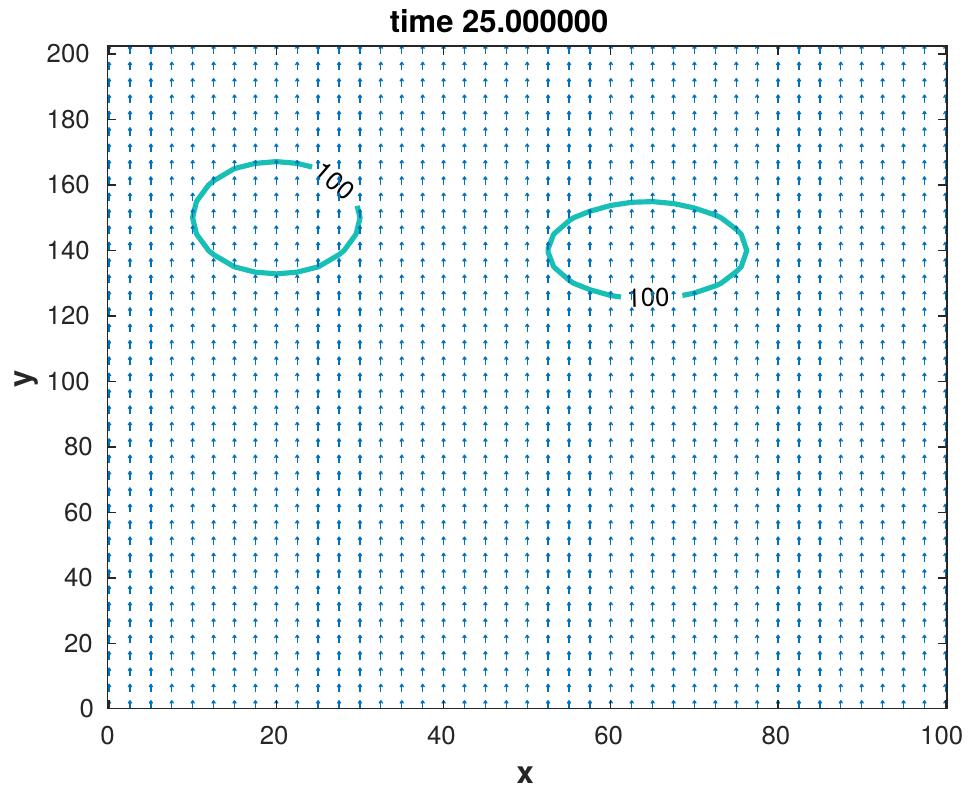}
\includegraphics[width=5.0cm,angle=-0]{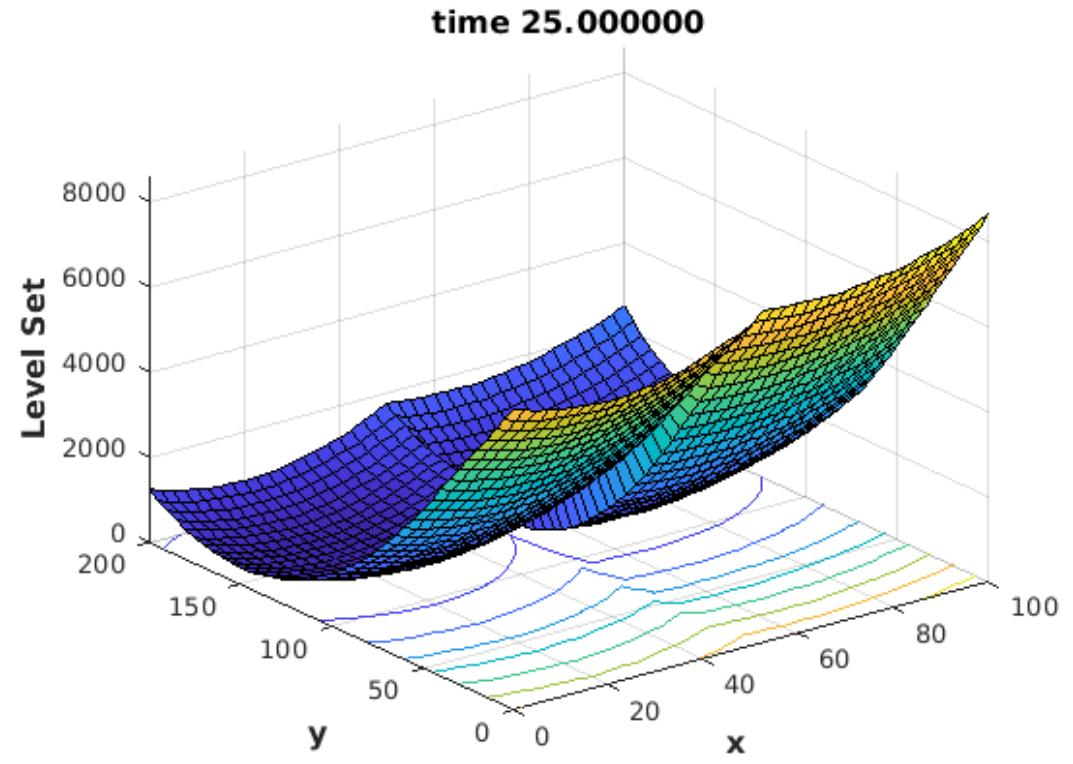}
\end{center}
\caption{\label{bubble_far_2bubb} Upper figures: Transport of the two bubble with the level-set function at time $t=0.1$, lower figures: Transport of the two bubble with the level-set function at time $t=25$.}
\end{figure}

The numerical results of the near-far-field coupled bubble-transport code,
which is given in Figure \ref{bubble_far_2bubb}.

\begin{remark}
The coupling of the formation and transport of the bubbles are done with ordinary and partial differential equations.
Based on decoupling such systems of mixed ordinary and partial differential equations, we
could compute each separate part with the optimal numerical solvers.
\end{remark}

\subsection{Bubble-Formation: Multiple Bubble Experiment (10 Bubbles)}
\label{bubble_10}

In the following, we extend the near-field and far-field computations
with 10 bubbles.
We also apply the decomposition of near-field and far-field computations 
as given in Figure \ref{near-far_1}.

We apply the following parameters:
\begin{itemize}
\item Computation of the near-field bubbles (a representing bubble is computed)
\begin{itemize}
\item Input-parameters of the near-field bubbles computation:
\begin{itemize}
\item Bubble 1: \\
 $r_1(0)=0.01, \; z_1(0)=0, \; r_1(L) =0, \; L = 2, \; \frac{\Delta p_1}{\alpha_1}= 0.2$, 
\item Bubble 2: \\
 $r_2(0)=0.01, \; z_2(0)=0, \; r_2(L) =0, \; L = 2, \; \frac{\Delta p_2}{\alpha_2}= 0.4$.
\item Bubble 3: \\
 $r_3(0)=0.01, \; z_3(0)=0, \; r_3(L) =0, \; L = 2, \; \frac{\Delta p_3}{\alpha_3}= 0.6$.
\item Bubble 4: \\
 $r_4(0)=0.01, \; z_4(0)=0, \; r_4(L) =0, \; L = 2, \; \frac{\Delta p_4}{\alpha_4}= 0.8$.
\item Bubble 5: \\
 $r_5(0)=0.01, \; z_5(0)=0, \; r_5(L) =0, \; L = 2, \; \frac{\Delta p_5}{\alpha_5}= 1.0$.
\item Bubble 6: \\
 $r_6(0)=0.01, \; z_6(0)=0, \; r_6(L) =0, \; L = 2, \; \frac{\Delta p_6}{\alpha_6}= 1.2$.
\item Bubble 7: \\
 $r_7(0)=0.01, \; z_7(0)=0, \; r_7(L) =0, \; L = 2, \; \frac{\Delta p_7}{\alpha_7}= 1.4$.
\item Bubble 8: \\
 $r_8(0)=0.01, \; z_8(0)=0, \; r_8(L) =0, \; L = 2, \; \frac{\Delta p_8}{\alpha_8}= 1.6$.
\item Bubble 9: \\
 $r_9(0)=0.01, \; z_9(0)=0, \; r_9(L) =0, \; L = 2, \; \frac{\Delta p_9}{\alpha_9}= 1.8$.
\item Bubble 10: \\
  $r_{10}(0)=0.01, \; z_{10}(0)=0, \; r_{10}(L) =0, \; L = 2, \; \frac{\Delta p_{10}}{\alpha_{10}}= 2.0$.
\end{itemize}
\item Output-parameters of the near-field bubble computation: 
\begin{itemize}
\item Bubble 1: \\
$a_{bubble_1} = 0.5141, \; b_{bubble_1} = 0.9926$.
\item Bubble 2: \\
$a_{bubble_{2}} = 0.5219, \; b_{bubble_{2}} = 0.9718$.
\item Bubble 3: \\
$a_{bubble_{3}} = 0.5295, \; b_{bubble_{3}} = 0.9506$.
\item Bubble 4: \\
$a_{bubble_{4}} = 0.5369, \; b_{bubble_{4}} = 0.9289$.
\item Bubble 5: \\
$a_{bubble_{5}} = 0.5443, \; b_{bubble_{5}} = 0.9067$.
\item Bubble 6: \\
$a_{bubble_{6}} = 0.5514, \; b_{bubble_{6}} = 0.8841$.
\item Bubble 7: \\
$a_{bubble_{7}} = 0.5584, \; b_{bubble_{7}} = 0.8612$.
\item Bubble 8: \\
$a_{bubble_{8}} = 0.5651, \; b_{bubble_{8}} = 0.8382$.
\item Bubble 9: \\
$a_{bubble_{9}} = 0.5717, \; b_{bubble_{9}} = 0.8154$.
\item Bubble 10: \\
$a_{bubble_{10}} = 0.5780, \; b_{bubble_{10}} = 0.7928$.
\end{itemize}
\item Ellipse: $(x-x_{bubble_i})^2 + ((y-y_{bubble_i})*a_{bubble}/b_{bubble})^2 - a_{bubble}^2 $, \\
where $(x_{bubble_i}, y_{bubble_i})$ is the origin of the $i$-th bubble.
\end{itemize}
\item Computation of the far-field bubbles (level-set initialisation):
\begin{itemize}
\item Parameterisation of the level-set initial-function, e.g., two bubbles:
\begin{eqnarray}
\phi_0(x,y) = \left\{
\begin{array}{cc}
(x - x_{bubble_1})^2 + ((y-y_{bubble_1}) \frac{a_{bubble_1}}{b_{bubble_1}})^2 - a_{bubble_1}^2 , & \\
 a_x \le x \le 50 , \;  a_y \le x \le b_y , & \\
(x - x_{bubble_2})^2 + ((y-y_{bubble_2}) \frac{a_{bubble_2}}{b_{bubble_2}})^2 - a_{bubble_2}^2 , & \\
 50 \le x \le b_x, \;  a_y \le x \le b_y , & \\
\end{array}
\right.
\end{eqnarray}
where $(x_{bubble_1}, y_{bubble_1}) = (20, 50)$, $(x_{bubble_2}, y_{bubble_2}) = (80, 50)$
with the coordinates of the grid $(a_x, a_y) = (0, 0)$ and $b_x, b_y) = (100, 200)$.
\end{itemize}

\end{itemize}

\begin{figure}[ht]
\begin{center}  
\includegraphics[width=5.0cm,angle=-0]{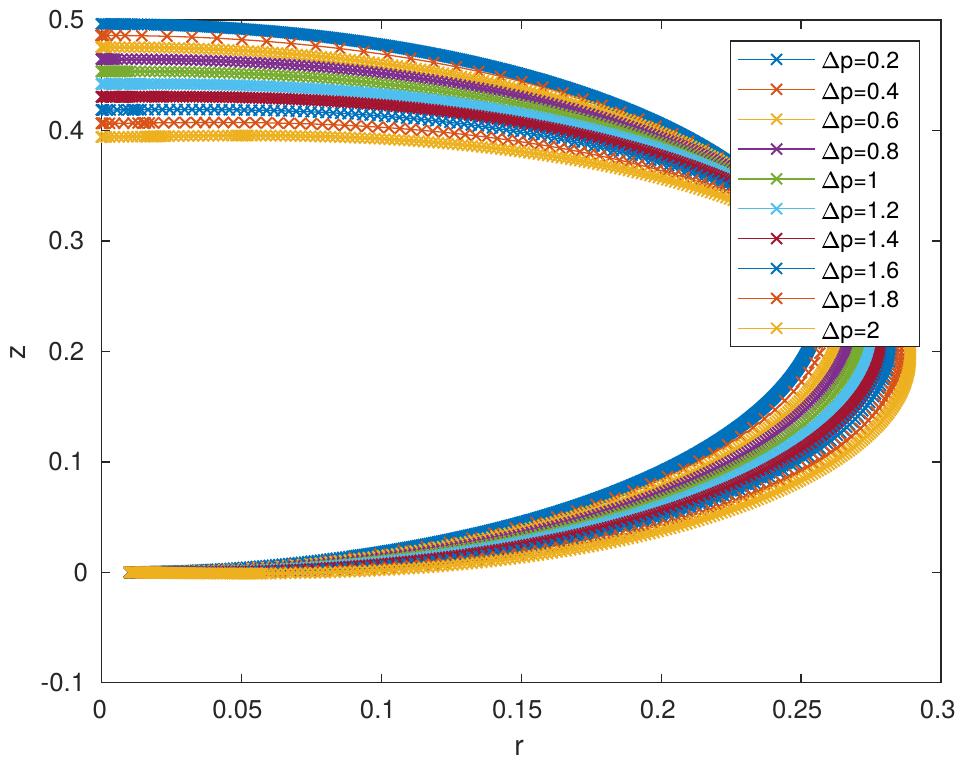}
\end{center}
\caption{\label{bubble_coupled_2_init} Formation of 10 bubbles with the different pressure-terms.}
\end{figure}

The numerical results of the formation of the bubbles are given in Figures \ref{bubble_coupled_2_init}.

\begin{figure}[ht]
\begin{center}  
\includegraphics[width=5.0cm,angle=-0]{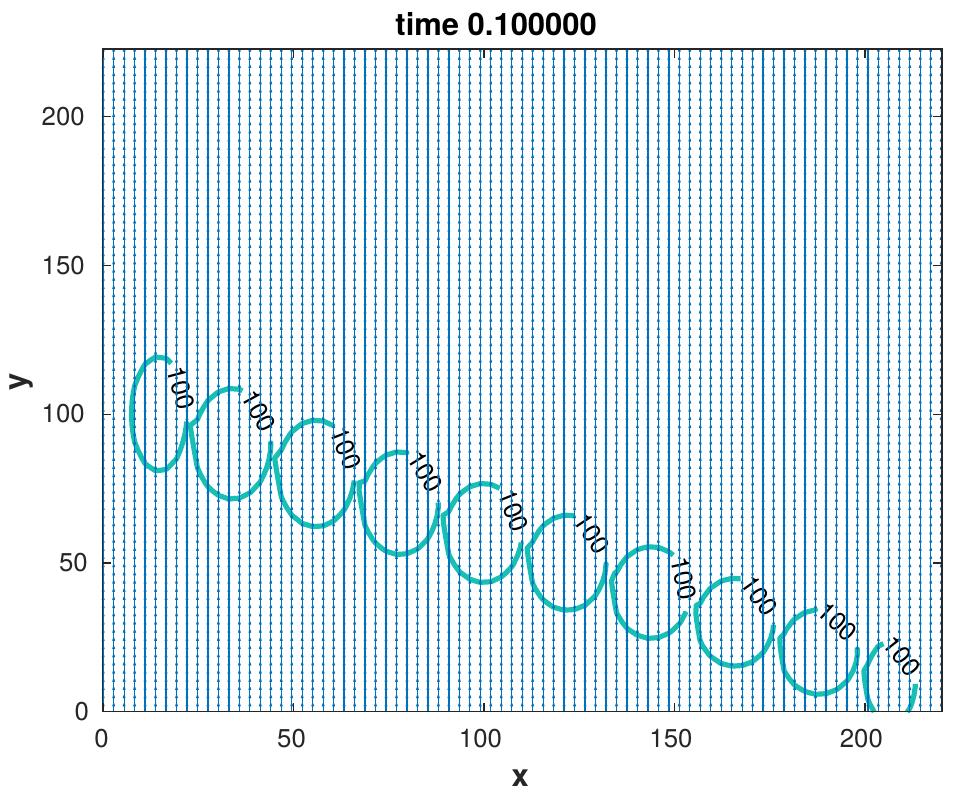}
\includegraphics[width=5.0cm,angle=-0]{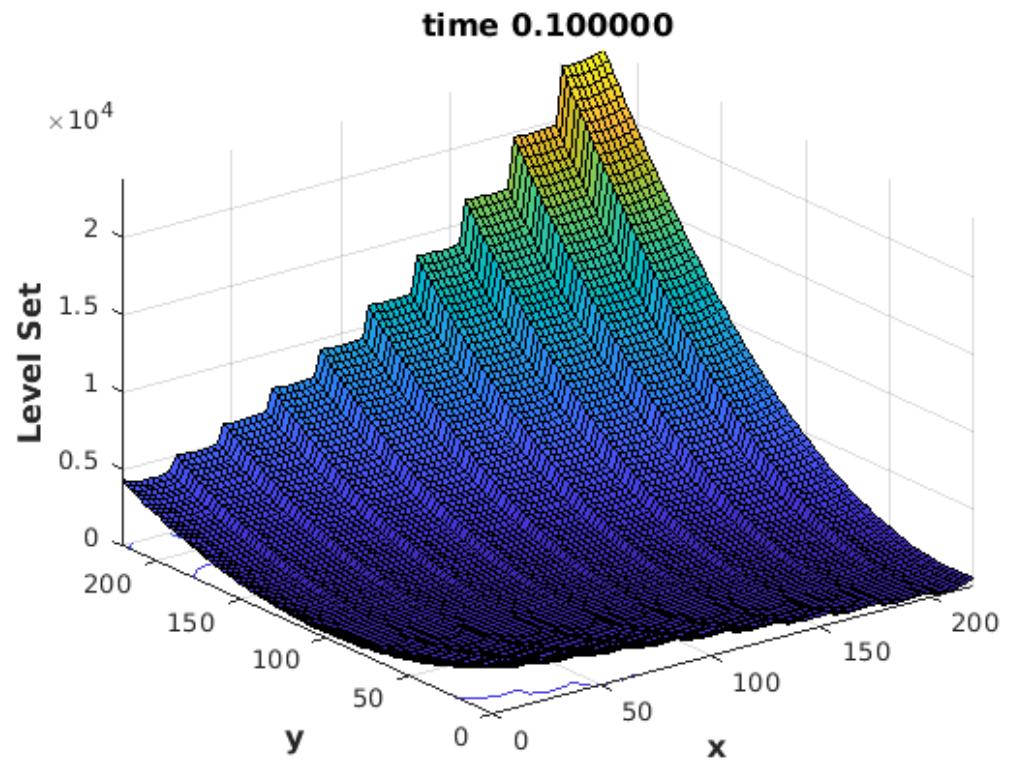}
\end{center}
\caption{\label{bubble_coupled_2} Transport of 10 bubbles with the level-set function at the initialisation $t = 0.1$.}
\end{figure}

\begin{figure}[ht]
\begin{center}  
\includegraphics[width=5.0cm,angle=-0]{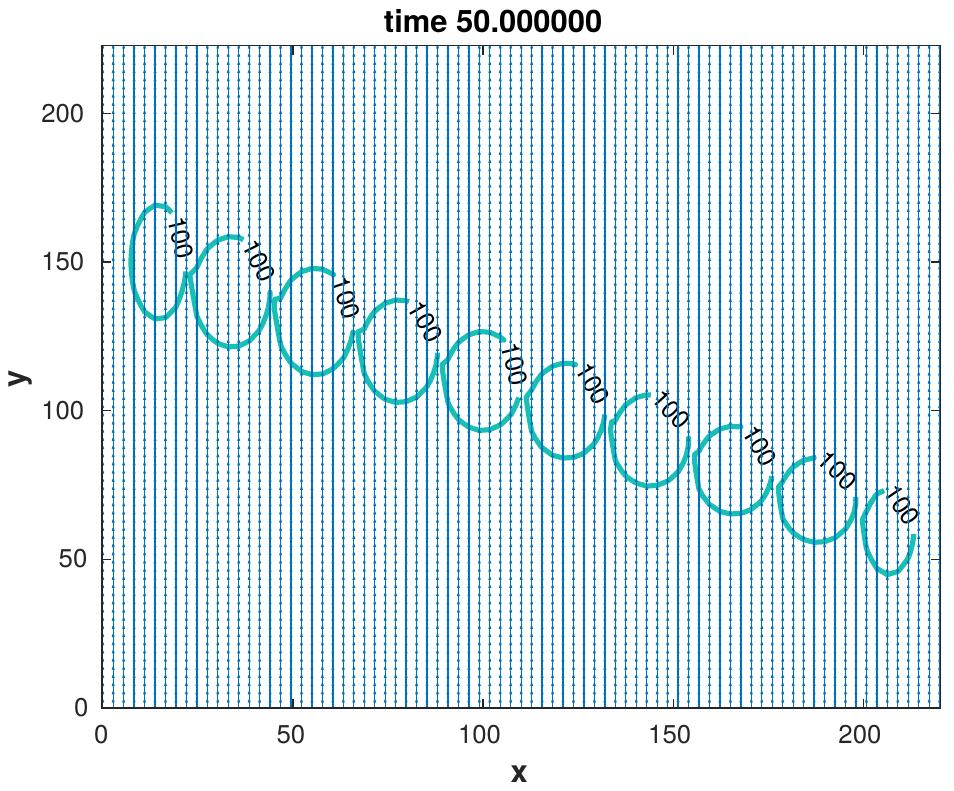}
\includegraphics[width=5.0cm,angle=-0]{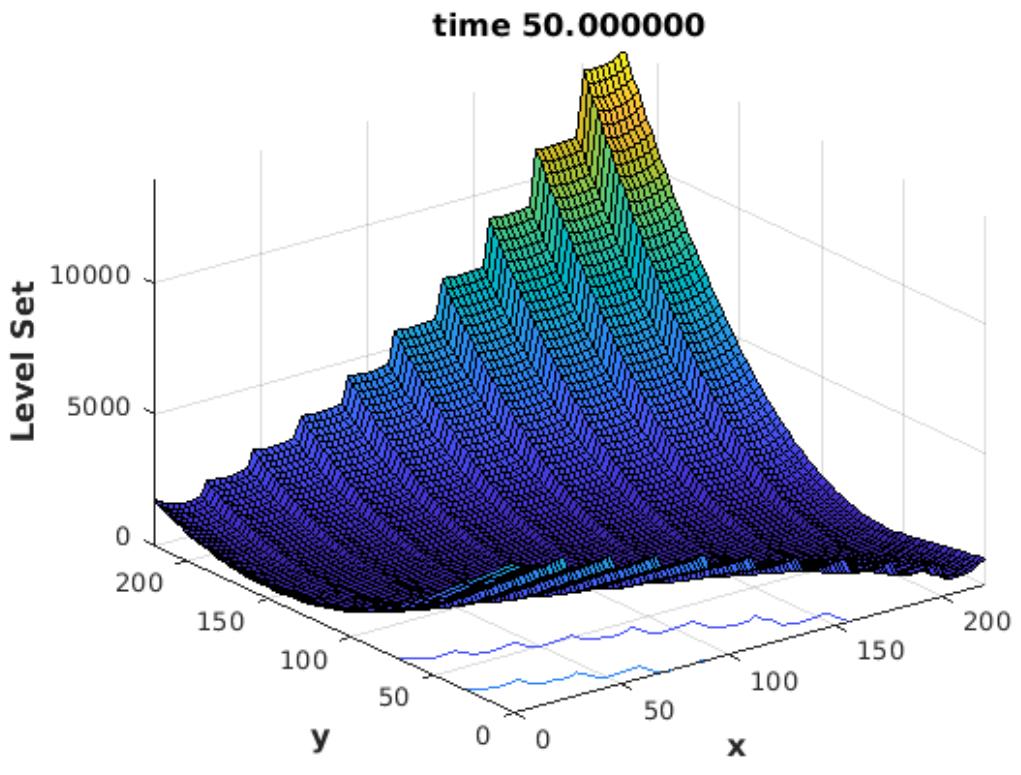}
\end{center}
\caption{\label{bubble_coupled_21} Transport of 10 bubbles with the level-set function at time $t=50$.}
\end{figure}

The numerical results of the near-far-field coupled bubble-transport code,
which is given in Figures \ref{bubble_coupled_2} and  \ref{bubble_coupled_21}.

\begin{remark}
In the experiment, we deal with at least 10 bubbles, which are 
different formated and transported via the level-set method.
The numerical experiments allows to accelerate the formation and transport
of such processes.
\end{remark}

\subsection{Bubble-Formation: Oscillation of the Air-Bubbles in the electrical-Field (first modelling approach)}

In the following, we simulate a bubble filled with air
in a electrical field, see \cite{sommers2012}.

\begin{figure}[ht]
\begin{center}  
\includegraphics[width=8.0cm,angle=-0]{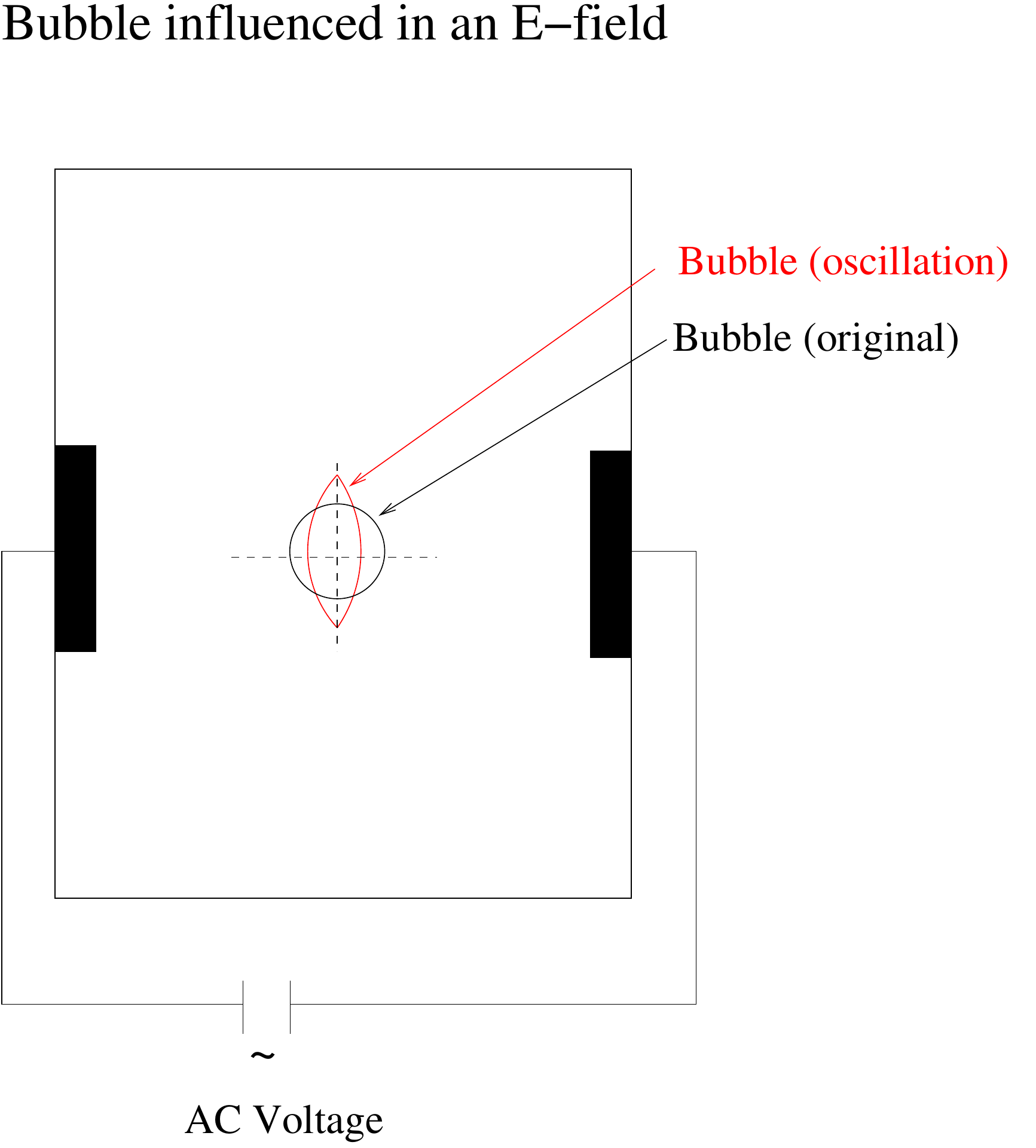}
\end{center}
\caption{\label{bubble_electr_1} Bubble-formation influenced by an electrical field.}
\end{figure}

We have the following setting of the influenced bubble, see Figure \ref{bubble_electr_1}.

We apply the following near-field equation based on an
extension of the pressure-term with an E-field.

The near-field equations are given as:
\begin{eqnarray}
  && \frac{d r}{d s} = \cos(\theta) , \\ 
  && \frac{d z}{d s} = \sin(\theta) , \\
  && \frac{d \theta}{d s} + \frac{\sin(\theta)}{r} = \frac{\Delta p}{\alpha} , 
\end{eqnarray}
where $s$ is the arc length along the curve
and $\theta$ the angle of elevation for its
slope and $\alpha$ is the mono-layer surface tension.

Further, we have 
\begin{eqnarray}
  && \Delta p = p_{bubble} - p_{tube} = \rho g z - p_{tube} , 
\end{eqnarray}
where $p_{tube} = p_E = \frac{9}{8} |{\bf E}_0|^2 \; \sin(\theta)$ is the pressure based on the electrical field of the tube
and we obtain:
\begin{eqnarray}
  && \frac{d r}{d s} = \cos(\theta) , \\ 
  && \frac{d z}{d s} = \sin(\theta) , \\
  && \frac{d \theta}{d s} + \frac{\sin(\theta)}{r} = \frac{\rho g z - \frac{9}{8} |{\bf E}_0|^2 \; \sin(\theta)}{\alpha} .
\end{eqnarray}

We apply the following parameters:
\begin{itemize}
\item Computation of the near-field bubbles (a representing bubble is computed)
\begin{itemize}
\item Input-parameters of the near-field bubbles computation:
\begin{itemize}
\item Bubble 1: \\
 $r_1(0)=0.01, \; z_1(0)=0, \; r_1(L) =0, \; L = 2, \; \frac{\Delta p_1}{\alpha_1}= 0.2$, 
\item Bubble 2: \\
 $r_2(0)=0.01, \; z_2(0)=0, \; r_2(L) =0, \; L = 2, \; \frac{\Delta p_2}{\alpha_2}= 0.4$.
\item Bubble 3: \\
 $r_3(0)=0.01, \; z_3(0)=0, \; r_3(L) =0, \; L = 2, \; \frac{\Delta p_3}{\alpha_3}= 0.6$.
\item Bubble 4: \\
 $r_4(0)=0.01, \; z_4(0)=0, \; r_4(L) =0, \; L = 2, \; \frac{\Delta p_4}{\alpha_4}= 0.8$.
\item Bubble 5: \\
 $r_5(0)=0.01, \; z_5(0)=0, \; r_5(L) =0, \; L = 2, \; \frac{\Delta p_5}{\alpha_5}= 1.0$.
\item Bubble 6: \\
 $r_6(0)=0.01, \; z_6(0)=0, \; r_6(L) =0, \; L = 2, \; \frac{\Delta p_6}{\alpha_6}= 1.2$.
\item Bubble 7: \\
 $r_7(0)=0.01, \; z_7(0)=0, \; r_7(L) =0, \; L = 2, \; \frac{\Delta p_7}{\alpha_7}= 1.4$.
\item Bubble 8: \\
 $r_8(0)=0.01, \; z_8(0)=0, \; r_8(L) =0, \; L = 2, \; \frac{\Delta p_8}{\alpha_8}= 1.6$.
\item Bubble 9: \\
 $r_9(0)=0.01, \; z_9(0)=0, \; r_9(L) =0, \; L = 2, \; \frac{\Delta p_9}{\alpha_9}= 1.8$.
\item Bubble 10: \\
  $r_{10}(0)=0.01, \; z_{10}(0)=0, \; r_{10}(L) =0, \; L = 2, \; \frac{\Delta p_{10}}{\alpha_{10}}= 2.0$.
\end{itemize}

\item Electrical field parameters are given as: 

$\alpha = 0.1$, $\rho = 0.1, \; g = 9.81, \;  |{\bf E}_0|^2 = 0.1$.

\item Output-parameters of the near-field bubble computation (formation): 
\begin{itemize}
\item Bubble 1: \\
$a_{bubble_1} = 0.5141, \; b_{bubble_1} = 0.9926$.
\item Bubble 2: \\
$a_{bubble_{2}} = 0.5219, \; b_{bubble_{2}} = 0.9718$.
\item Bubble 3: \\
$a_{bubble_{3}} = 0.5295, \; b_{bubble_{3}} = 0.9506$.
\item Bubble 4: \\
$a_{bubble_{4}} = 0.5369, \; b_{bubble_{4}} = 0.9289$.
\item Bubble 5: \\
$a_{bubble_{5}} = 0.5443, \; b_{bubble_{5}} = 0.9067$.
\item Bubble 6: \\
$a_{bubble_{6}} = 0.5514, \; b_{bubble_{6}} = 0.8841$.
\item Bubble 7: \\
$a_{bubble_{7}} = 0.5584, \; b_{bubble_{7}} = 0.8612$.
\item Bubble 8: \\
$a_{bubble_{8}} = 0.5651, \; b_{bubble_{8}} = 0.8382$.
\item Bubble 9: \\
$a_{bubble_{9}} = 0.5717, \; b_{bubble_{9}} = 0.8154$.
\item Bubble 10: \\
$a_{bubble_{10}} = 0.5780, \; b_{bubble_{10}} = 0.7928$.
\end{itemize}

\item Output-parameters of the near-field bubble computation (in the E-field): 
\begin{itemize}
\item Bubble 1: \\
$a_{bubble_1} = 0.5558, \; b_{bubble_1} = 0.8698$.
\item Bubble 2: \\
$a_{bubble_{2}} = 0.5626, \; b_{bubble_{2}} = 0.8468$.
\item Bubble 3: \\
$a_{bubble_{3}} = 0.5693, \; b_{bubble_{3}} = 0.8239$.
\item Bubble 4: \\
$a_{bubble_{4}} = 0.5757, \; b_{bubble_{4}} = 0.8013$.
\item Bubble 5: \\
$a_{bubble_{5}} = 0.5819, \; b_{bubble_{5}} = 0.7789$.
\item Bubble 6: \\
$a_{bubble_{6}} = 0.5878, \; b_{bubble_{6}} = 0.7571$.
\item Bubble 7: \\
$a_{bubble_{7}} = 0.5935, \; b_{bubble_{7}} = 0.7359$.
\item Bubble 8: \\
$a_{bubble_{8}} = 0.5990, \; b_{bubble_{8}} = 0.7153$.
\item Bubble 9: \\
$a_{bubble_{9}} = 0.6042, \; b_{bubble_{9}} = 0.6953$.
\item Bubble 10: \\
$a_{bubble_{10}} = 0.6091, \; b_{bubble_{10}} = 0.6760$.
\end{itemize}

\item Ellipse: $(x-x_{bubble_i})^2 + ((y-y_{bubble_i})*a_{bubble}/b_{bubble})^2 - a_{bubble}^2 $, \\
where $(x_{bubble_i}, y_{bubble_i})$ is the origin of the $i$-th bubble.
\end{itemize}
\item Computation of the far-field bubbles (level-set initialisation):
\begin{itemize}
\item Parameterisation of the level-set initial-function, e.g., two bubbles:
\begin{eqnarray}
\phi_0(x,y) = \left\{
\begin{array}{cc}
(x - x_{bubble_1})^2 + ((y-y_{bubble_1}) \frac{a_{bubble_1}}{b_{bubble_1}})^2 - a_{bubble_1}^2 , & \\
 a_x \le x \le 50 , \;  a_y \le x \le b_y , & \\
(x - x_{bubble_2})^2 + ((y-y_{bubble_2}) \frac{a_{bubble_2}}{b_{bubble_2}})^2 - a_{bubble_2}^2 , & \\
 50 \le x \le b_x, \;  a_y \le x \le b_y , & \\
\end{array}
\right.
\end{eqnarray}
where $(x_{bubble_1}, y_{bubble_1}) = (20, 50)$, $(x_{bubble_2}, y_{bubble_2}) = (80, 50)$
with the coordinates of the grid $(a_x, a_y) = (0, 0)$ and $b_x, b_y) = (100, 200)$.
\end{itemize}

\end{itemize}

\begin{figure}[ht]
\begin{center}  
\includegraphics[width=5.0cm,angle=-0]{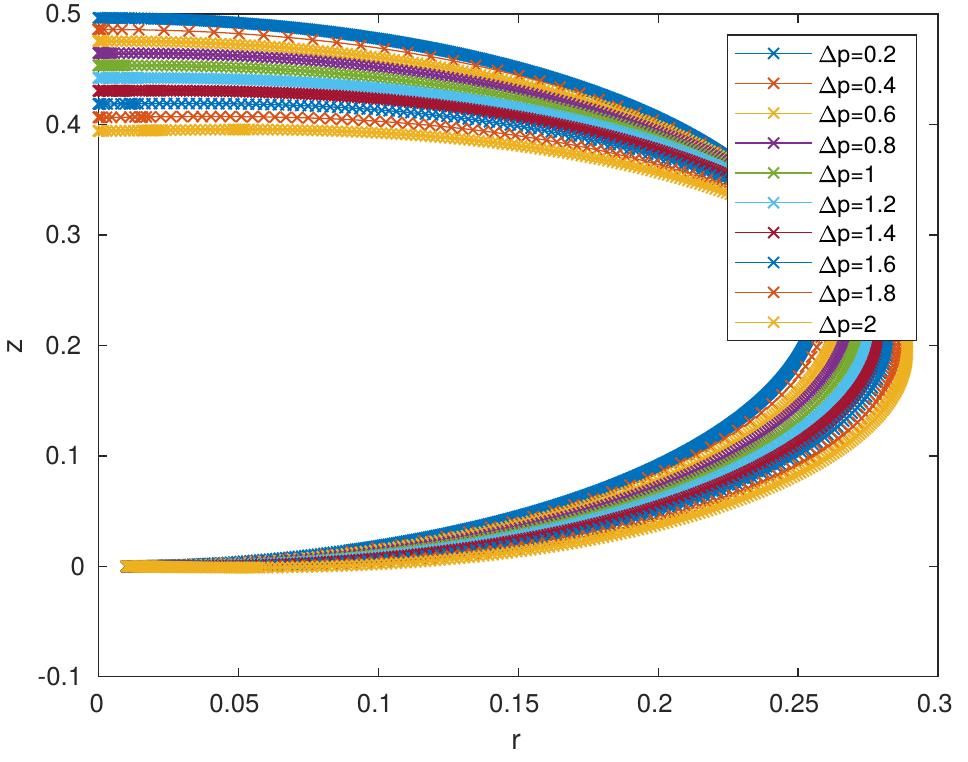}
\includegraphics[width=5.0cm,angle=-0]{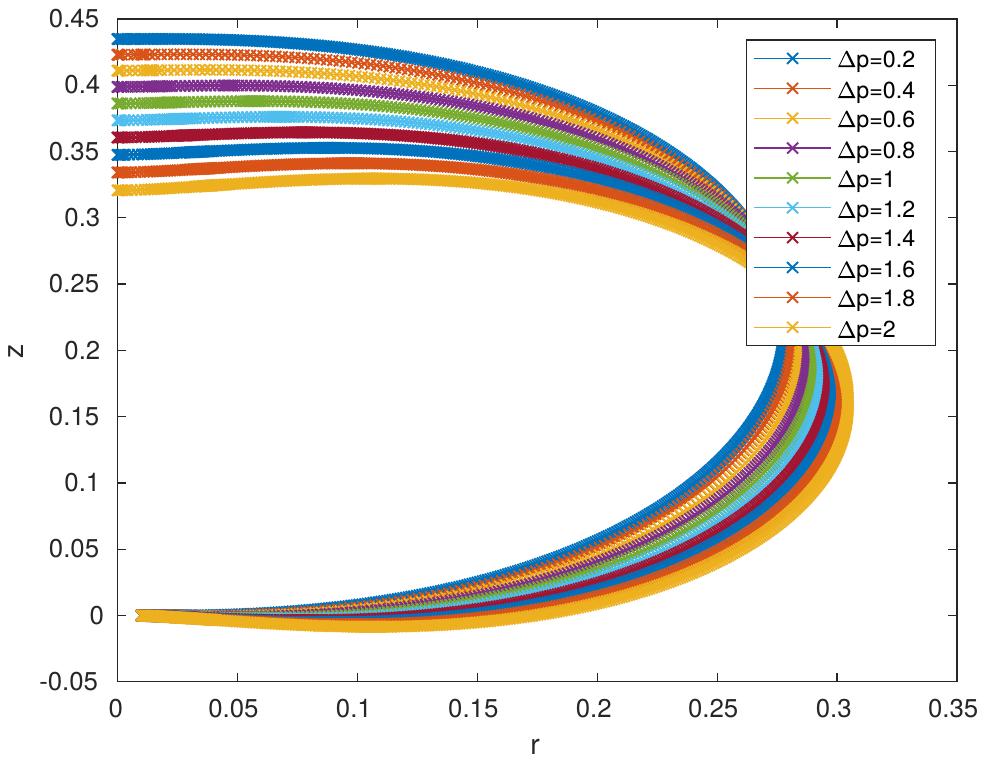}
\end{center}
\caption{\label{bubble_el_0_init} Left figure: Formation of 10 bubbles without an E-field, right figure: Formation of 10 bubbles with an E-field.}
\end{figure}

The numerical results of the formation of the bubbles with and without the E-field are given in Figures \ref{bubble_el_0_init}.

\begin{figure}[ht]
\begin{center}  
\includegraphics[width=5.0cm,angle=-0]{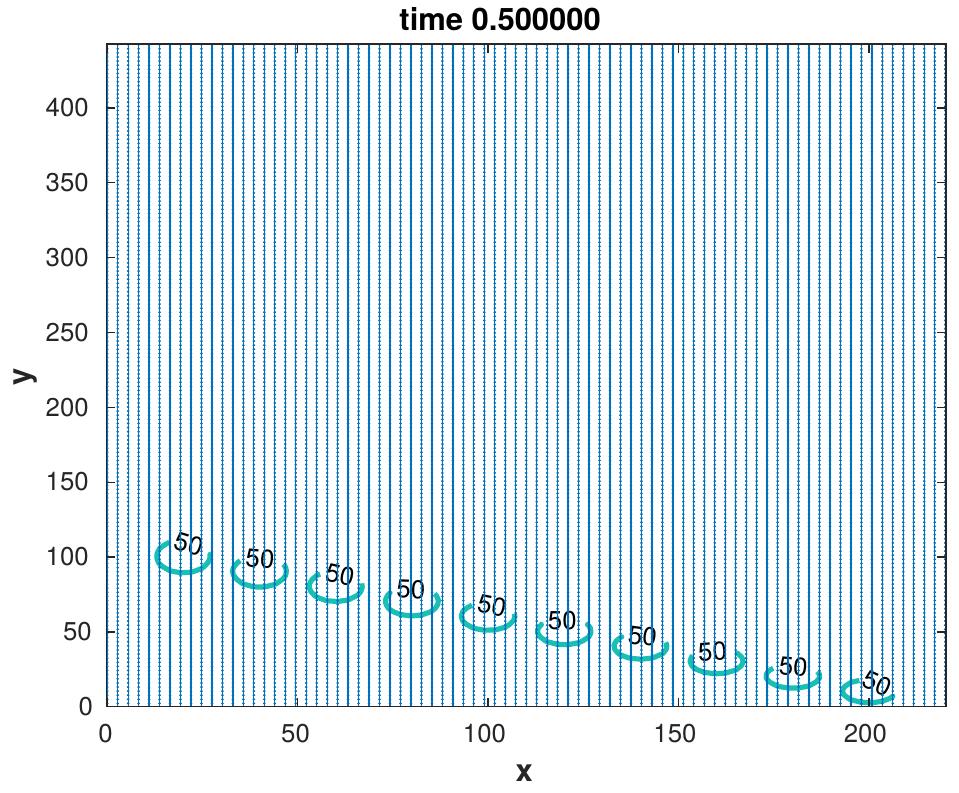}
\includegraphics[width=5.0cm,angle=-0]{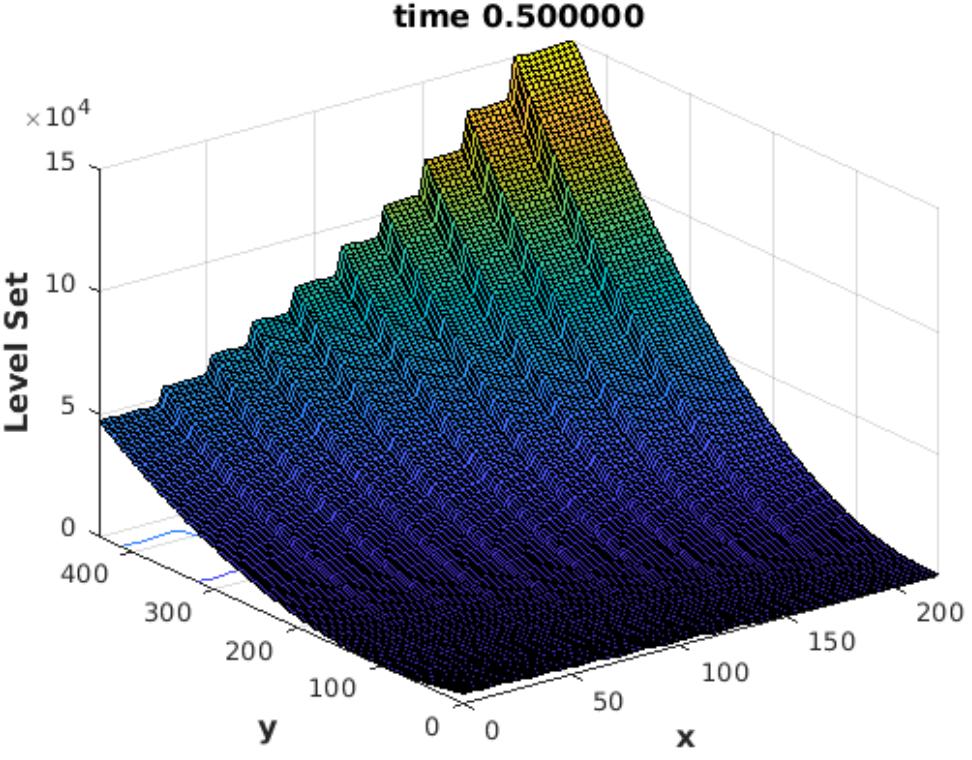}
\end{center}
\caption{\label{bubble_coupled_elect_2} Transport of 10 bubbles in the E-field with the level-set function at the initialisation $t = 0.1$.}
\end{figure}

\begin{figure}[ht]
\begin{center}  
\includegraphics[width=5.0cm,angle=-0]{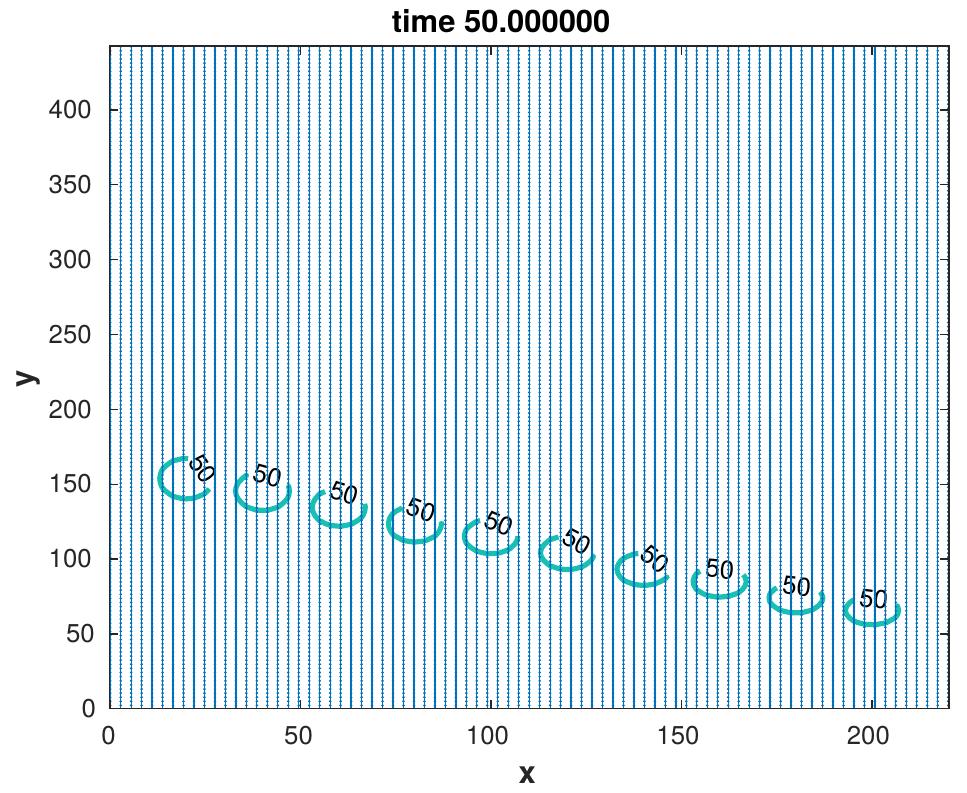}
\includegraphics[width=5.0cm,angle=-0]{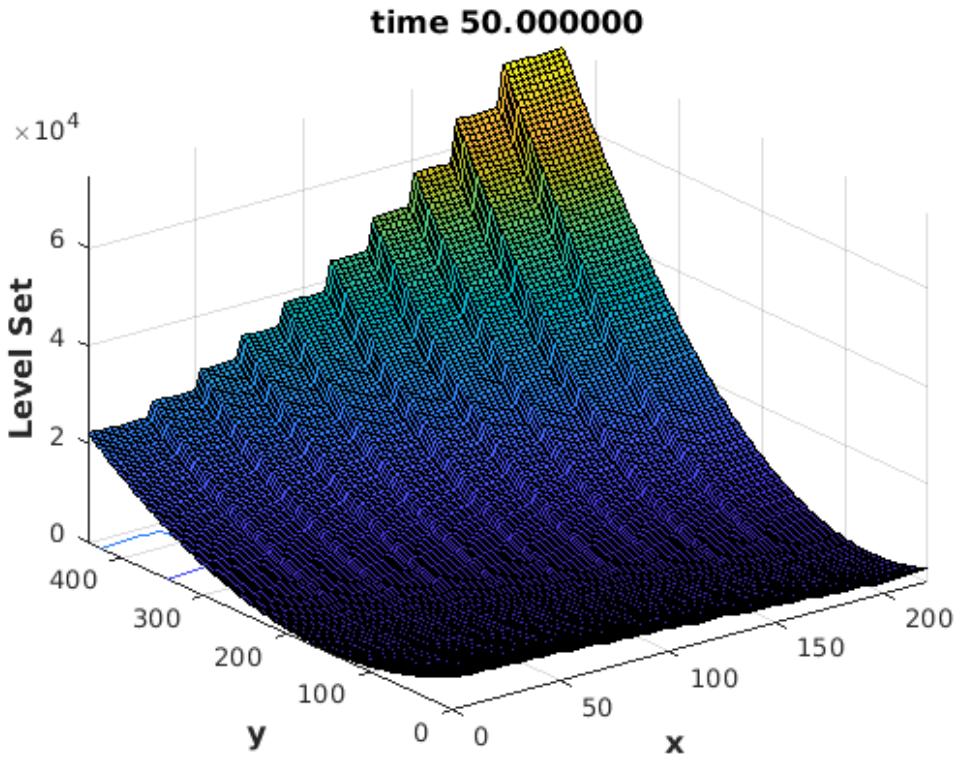}
\end{center}
\caption{\label{bubble_coupled_elect_21} Transport of 10 bubbles in the E-field with the level-set function at time $t=50$.}
\end{figure}

The numerical results of the near-far-field coupled bubble-transport code,
which is given in Figures \ref{bubble_coupled_elect_2} and  \ref{bubble_coupled_elect_21}.

\begin{remark}
The bubble modifications are given by the E-field, while it changes the formation of the bubble.
The level-set method is modified by multi-level-set domains and allows to deal with different level-set functions, such that we could 
transport multiple bubbles with E-fields.
In the first approach, we assume, that the E-field can be approximated at the beginning of the formation and that the bubble will not be changed after such a initial formation.
\end{remark}

\section{Conclusion}
\label{concl}

We present a bubble model, which is a coupled model based on a bubble formation and bubble transport
model. The decoupling into near- and far-field models allows to apply optimal solver and
discretization methods. We apply different numerical experiments, which shows the benefit of
such a treatment. In future, we consider the fully coupled problem, while we deal with
bubble density functions and the coupling between the formation and transport process.
Such an extension allows to see the ruptures of the bubbles.

\bibliographystyle{plain}

\end{document}